# FUNCTIONAL LINEAR REGRESSION ANALYSIS FOR LONGITUDINAL DATA[1]

By Fang Yao, Hans-Georg Müller and Jane-Ling Wang

*Colorado State University, University of California, Davis, and University of California, Davis*

We propose nonparametric methods for functional linear regression which are designed for sparse longitudinal data, where both the predictor and response are functions of a covariate such as time. Predictor and response processes have smooth random trajectories, and the data consist of a small number of noisy repeated measurements made at irregular times for a sample of subjects. In longitudinal studies, the number of repeated measurements per subject is often small and may be modeled as a discrete random number and, accordingly, only a finite and asymptotically nonincreasing number of measurements are available for each subject or experimental unit. We propose a functional regression approach for this situation, using functional principal component analysis, where we estimate the functional principal component scores through conditional expectations. This allows the prediction of an unobserved response trajectory from sparse measurements of a predictor trajectory. The resulting technique is flexible and allows for different patterns regarding the timing of the measurements obtained for predictor and response trajectories. Asymptotic properties for a sample of $n$ subjects are investigated under mild conditions, as $n \to \infty$, and we obtain consistent estimation for the regression function. Besides convergence results for the components of functional linear regression, such as the regression parameter function, we construct asymptotic pointwise confidence bands for the predicted trajectories. A functional coefficient of determination as a measure of the variance explained by the functional regression model is introduced, extending the standard $R^2$ to the functional case. The proposed methods are illustrated with a simulation study, longitudinal primary biliary liver cirrhosis data and an analysis of

Received June 2004; revised December 2004.

[1]Supported in part by NSF Grants DMS-99-71602, DMS-00-79430, DMS-02-04869, DMS-03-54448 and DMS-04-06430.

the longitudinal relationship between blood pressure and body mass index.

**1. Introduction.** We develop a version of functional linear regression analysis in which both the predictor and response variables are functions of some covariate which usually but not necessarily is time. Our approach extends the applicability of functional regression to typical longitudinal data where only very few and irregularly spaced measurements for predictor and response functions are available for most of the subjects. Examples of such data are discussed in Section 5 (see Figures 1 and 6).

Since a parametric approach only captures features contained in the pre-conceived class of functions, nonparametric methods of functional data analysis are needed for the detection of new features and for the modeling of highly complex relationships. Functional principal component analysis (FPCA) is a basic methodology that has been studied in early work by Grenander [18] and, more recently, by Rice and Silverman [27], Ramsay and Silverman [26] and many others. Background in probability on function spaces can be found in [19]. James, Hastie and Sugar [21] emphasized the case of sparse data by proposing a reduced rank mixed-effects model using B-spline functions. Nonparametric methods for unbalanced longitudinal data were studied by Boularan, Ferré and Vieu [2] and Besse, Cardot and Ferraty [1]. Yao, Müller and Wang [31] proposed an FPCA procedure through a conditional expectation method, aiming at estimating functional principal component scores for sparse longitudinal data.

In the recent literature there has been increased interest in regression models for functional data, where both the predictor and response are random functions. Our aim is to extend the applicability of such models to longitudinal data with their typically irregular designs, and to develop asymptotics for functional regression in sparse data situations. Practically all investigations to date are for the case of completely observed trajectories, where one assumes either entire trajectories or densely spaced measurements taken along each trajectory are observed; recent work includes Cardot, Ferraty, Mas and Sarda [3], Cardot, Ferraty and Sarda [5], Chiou, Müller, Wang and Carey [7] and Ferraty and Vieu [16].

In this paper we illustrate the potential of functional regression for complex longitudinal data. In functional data settings, Cardot, Ferraty and Sarda [4] provided consistency results for the case of linear regression with a functional predictor and scalar response, where the predictor functions are sampled at a regular grid for each subject, and Cardot, Ferraty and Sarda [4] discussed inference for the regression function. The case of a functional response was introduced by Ramsay and Dalzell [25], and for a summary of this and related work we refer to Ramsay and Silverman ([26], Chapter 11) and to Faraway [15] for a discussion of relevant practical aspects. The theory for the case of fixed design and functional response in the densely sampled case was investigated by Cuevas, Febrero and Fraiman [8]. Chiou, Müller and Wang [6] studied functional regression models where the predictors are finite-dimensional vectors and the response is a function, using a quasi-likelihood approach. Applications of varying-coefficient modeling to functional data, including asymptotic inference, were presented in [13] and [14].

The proposed functional regression approach is flexible, and allows for varying patterns of timing in regard to the measurements of predictor and response functions. This is relevant since it is a common occurrence in longitudinal data settings that the measurement of either the predictor or response is missing. The contributions of this paper are as follows: First, we extend the functional regression approach to longitudinal data, using a conditioning idea. This leads to improved prediction of the response trajectories, given sparse measurements of the predictor trajectories. Second, we provide a complete practical implementation of the proposed functional regression procedure and illustrate its utility for two longitudinal studies. Third, we obtain the asymptotic consistency of the estimated regression function of the functional linear regression model for the case of sparse and irregular data, including rates. Fourth, we construct asymptotic pointwise confidence bands for predicted response trajectories based on asymptotic distribution results. Fifth, we introduce a consistent estimator for a proposed measure of association between the predictor and response functions in functional regression models that provides an extension of the coefficient of determination $R^2$ in standard linear model theory to the functional case. The proposed functional coefficient of determination provides a useful quantification of the strength of the relationship between response and predictor functions, as it can be interpreted in a well-defined sense as the fraction of variance explained by the functional linear regression model, in analogy to the situation for the standard linear regression model.

The paper is organized as follows. In Section 2 we introduce basic notions, the functional linear regression model, and describe the estimation of the regression function. In Section 3 we discuss the extension of the conditioning approach to the prediction of response trajectories in functional regression under irregular and sparse data. Pointwise confidence bands and the functional coefficient of determination $R^2$ are also presented in Section 3. Simulation results that illustrate the usefulness of the proposed method can be found in Section 4. This is followed by applications of the proposed functional regression approach to longitudinal PBC liver cirrhosis data and an analysis of the longitudinal relationship between blood pressure and body mass index, using data from the Baltimore Longitudinal Study on Aging in Section 5. Asymptotic consistency and distribution results are provided in Section 6, while proofs and auxiliary results are compiled in the Appendix.

## 2. Functional linear regression for sparse and irregular data.

2.1. *Representing predictor and response functions through functional principal components.* The underlying but unobservable sample consists of pairs of random trajectories $(X_i, Y_i)$, $i = 1, \ldots, n$, with square integrable predictor trajectories $X_i$ and response trajectories $Y_i$. These are realizations of smooth random processes $(X, Y)$, with unknown smooth mean functions $EY(t) = \mu_Y(t)$, $EX(s) = \mu_X(s)$, and covariance functions $\text{cov}(Y(s), Y(t)) = G_Y(s, t)$, $\text{cov}(X(s), X(t)) = G_X(s, t)$. We usually refer to the arguments of $X(\cdot)$ and $Y(\cdot)$ as time, with finite and closed intervals $\mathcal{S}$ and $\mathcal{T}$ as domains. We assume the existence of orthogonal expansions of $G_X$ and $G_Y$ (in the $L^2$ sense) in terms of eigenfunctions $\psi_m$ and $\phi_k$ with nonincreasing eigenvalues $\rho_m$ and $\lambda_k$, that is, $G_X(s_1, s_2) = \sum \rho_m \psi_m(s_1) \psi_m(s_2)$, $t, s \in \mathcal{S}$, and $G_Y(t_1, t_2) = \sum_k \lambda_k \phi_k(t_1) \phi_k(t_2)$, $t_1, t_2 \in \mathcal{T}$.

We model the actually observed data which consist of sparse and irregular repeated measurements of the predictor and response trajectories $X_i$ and $Y_i$, contaminated with additional measurement errors (see [28, 30]). To adequately reflect the irregular and sparse measurements, we assume that there is a random number of $L_i$ (resp. $N_i$) random measurement times for $X_i$ (resp. $Y_i$) for the $i$th subject, which are denoted as $S_{i1}, \ldots, S_{iL_i}$ (resp. $T_{i1}, \ldots, T_{iN_i}$). The random variables $L_i$ and $N_i$ are assumed to be i.i.d. as $L$ and $N$ respectively, where $L$ and $N$ may be correlated but are independent of all other random variables. Let $U_{il}$ (resp. $V_{ij}$) denote the observation of the random trajectory $X_i$ (resp. $Y_i$) at a random time $S_{il}$ (resp. $T_{ij}$), contaminated with measurement errors $\varepsilon_{il}$ (resp. $\epsilon_{ij}$), $1 \le l \le L_i$, $1 \le j \le N_i$, $1 \le i \le n$. The errors are assumed to be i.i.d. with $E\varepsilon_{il} = 0$, $E[\varepsilon_{il}^2] = \sigma_X^2$ (resp. $E\epsilon_{ij} = 0$, $E[\epsilon_{ij}^2] = \sigma_Y^2$), and independent of functional principal component scores $\zeta_{im}$ (resp. $\xi_{ik}$) that satisfy $E\zeta_{im} = 0$, $E[\zeta_{im}\zeta_{im'}] = 0$ for $m \ne m'$, $E[\zeta_{im}^2] = \rho_m$ (resp. $E\xi_{ik} = 0$, $E[\xi_{ik}\xi_{ik'}] = 0$ for $k \ne k'$, $E[\xi_{ik}^2] = \lambda_k$). Then we may represent predictor and response measurements as follows:

$$\begin{aligned}
(1) \quad U_{il} &= X_i(S_{il}) + \varepsilon_{il} \\
&= \mu_X(S_{il}) + \sum_{m=1}^{\infty} \zeta_{im} \psi_m(S_{il}) + \varepsilon_{il}, \qquad S_{il} \in \mathcal{S}, 1 \le i \le n, 1 \le l \le L_i,
\end{aligned}$$

$$\begin{aligned}
(2) \quad V_{ij} &= Y_i(T_{ij}) + \epsilon_{ij} \\
&= \mu_Y(T_{ij}) + \sum_{k=1}^{\infty} \xi_{ik} \phi_k(T_{ij}) + \epsilon_{ij}, \qquad T_{ij} \in \mathcal{T}, 1 \le i \le n, 1 \le j \le N_i.
\end{aligned}$$

We note that the response and predictor functions do not need to be sampled simultaneously, extending the applicability of the proposed functional regression model.

2.2. *Functional linear regression model and estimation of the regression function.* Consider a functional linear regression model in which both the predictor $X$ and response $Y$ are smooth random functions,

$$(3) \qquad E[Y(t)|X] = \alpha(t) + \int_{\mathcal{S}} \beta(s,t)X(s)\,ds.$$

Here the bivariate regression function $\beta(s,t)$ is smooth and square integrable, that is, $\int_{\mathcal{T}} \int_{\mathcal{S}} \beta^2(s,t)\,ds\,dt < \infty$. Centralizing $X(t)$ by $X^c(s) = X(s) - \mu_X(s)$, and observing $E[Y(t)] = \mu_Y(t) = \alpha(t) + \int_{\mathcal{S}} \beta(s,t)\mu_X(s)\,ds$, the functional linear regression model becomes

$$(4) \qquad E[Y(t)|X] = \mu_Y(t) + \int_{\mathcal{S}} \beta(s,t)X^c(s)\,ds.$$

Our aim is to predict an unknown response trajectory based on sparse and noisy observations of a new predictor function. This is the functional version of the classical prediction problem in a linear model where, given a set of predictors $X$, one aims to predict the mean response $Y$ by estimating $E(Y|X)$ (see [12], page 81). An important step is to estimate the regression function $\beta(s,t)$. We use the following basis representation of $\beta(s,t)$, which is a consequence of the population least squares property of conditional expectation and the fact that the predictors are uncorrelated, generalizing the representation $\beta_1 = \operatorname{cov}(X,Y)/\operatorname{var}(X)$ of the slope parameter in the simple linear regression model $E(Y|X) = \beta_0 + \beta_1 X$ to the functional case. This representation holds under certain regularity conditions which are outlined in [20] and is given by

$$(5) \qquad \beta(s,t) = \sum_{k=1}^{\infty} \sum_{m=1}^{\infty} \frac{E[\zeta_m \xi_k]}{E[\zeta_m^2]} \psi_m(s)\phi_k(t).$$

The convergence of the right-hand side of (5) is discussed in Lemma A.2 (Appendix A.3). When referring to $\beta$, we always assume that the limit (5) exists in an appropriate sense. In a first step, smooth estimates of the mean and covariance functions for the predictor and response functions are obtained by scatterplot smoothing; see (30) and (31) in Appendix A.2. Then a nonparametric FPCA step yields estimates $\hat{\psi}_m$, $\hat{\phi}_k$ for the eigenfunctions, and $\hat{\rho}_m$, $\hat{\lambda}_k$ for the eigenvalues of predictor and response functions; see (33) below.

We use two-dimensional scatterplot smoothing to obtain an estimate $\widehat{C}(s,t)$ of the cross-covariance surface $C(s,t)$, $s \in \mathcal{S}$, $t \in \mathcal{T}$,

$$(6) \qquad C(s,t) = \operatorname{cov}(X(s),Y(t)) = \sum_{k=1}^{\infty} \sum_{m=1}^{\infty} E[\zeta_m \xi_k]\psi_m(s)\phi_k(t).$$

Let $C_i(S_{il}, T_{ij}) = (U_{il} - \hat{\mu}_X(S_{il}))(V_{ij} - \hat{\mu}_Y(T_{ij}))$ be "raw" cross-covariances that serve as input for the two-dimensional smoothing step; see (36) in Appendix A.2. The smoothing parameters in the two coordinate directions can be chosen independently by one-curve-leave-out cross-validation procedures [27]. From (6) we obtain estimates for $\sigma_{km} = E[\zeta_m \xi_k]$,

$$\hat{\sigma}_{km} = \int_{\mathcal{T}} \int_{\mathcal{S}} \hat{\psi}_m(s) \widehat{C}(s,t) \hat{\phi}_k(t) \, ds \, dt,$$

$$m = 1, \ldots, M, k = 1, \ldots, K.$$

(7)

With estimates (33), the resulting estimate for $\beta(s,t)$ is

$$\hat{\beta}(s,t) = \sum_{k=1}^{K} \sum_{m=1}^{M} \frac{\hat{\sigma}_{km}}{\hat{\rho}_m} \hat{\psi}_m(s) \hat{\phi}_k(t).$$

(8)

In practice, the numbers $M$ and $K$ of included eigenfunctions can be chosen by one-curve-leave-out cross-validation (34), or by an AIC type criterion (35). For the asymptotic analysis, we consider $M(n), K(n) \to \infty$ as the sample size $n \to \infty$. Corresponding convergence results can be found in Theorem 1.

## 3. Prediction and inference.

3.1. *Predicting response trajectories.* One of our central aims is to predict the trajectory $Y^*$ of the response for a new subject from sparse and irregular measurements of the predictor trajectory $X^*$. In view of (4), the basis representation of $\beta(s,t)$ in (5) and the orthonormality of the $\{\psi_m\}_{m \geq 1}$, the prediction of the response function would be obtained via the conditional expectation

$$E[Y^*(t)|X^*] = \mu_Y(t) + \sum_{k=1}^{\infty} \sum_{m=1}^{\infty} \frac{\sigma_{km}}{\rho_m} \zeta_m^* \phi_k(t),$$

(9)

where $\zeta_m^* = \int_{\mathcal{S}} (X^*(s) - \mu_X(s)) \psi_m(s) \, ds$ is the $m$th functional principal component score of the predictor trajectory $X^*$. The quantities $\mu_Y$, $\phi_k$, $\sigma_{km}$ and $\rho_m$ can be estimated from the data, as described above. It remains to discuss the estimation of $\zeta_m^*$, and for this step we invoke Gaussian assumptions in order to handle the sparsity of the data.

Let $U_l^*$ be the $l$th measurement made for the predictor function $X^*$ at time $S_l^*$, according to (1), where $l = 1, \ldots, L^*$, with $L^*$ a random number. Assume that the functional principal component scores $\zeta_m^*$ and the measurement errors $\varepsilon_l^*$ for the predictor trajectories are jointly Gaussian. Following Yao, Müller and Wang [31], the best prediction of the scores $\zeta_m^*$ is then obtained through the best linear prediction, given the observations $\widetilde{U}^* = (U_1^*, \ldots, U_{L^*}^*)$, and the number and locations of these observations, $L^*$ and $S^* = (S_1^*, \ldots, S_{L^*}^*)^T$. Let $X^*(S_l^*)$ be the value of the predictor function $X^*$ at time $S_l^*$. Write $\widetilde{X}^* = (X^*(S_1^*), \ldots, X^*(S_{L^*}^*))^T$, $\mu_X^* = (\mu_X(S_1^*), \ldots, \mu_X(S_{L^*}^*))^T$ and $\psi_m^* = (\psi_m(S_1^*), \ldots, \psi_m(S_{L^*}^*))^T$. Then the best linear prediction for $\zeta_m^*$ is

$$\tilde{\zeta}_m^* = \rho_m \psi_m^{*T} \Sigma_{U^*}^{-1} (\tilde{U}^* - \mu_X^*), \tag{10}$$

where $\Sigma_{U^*} = \text{cov}(\tilde{U}^*|L^*, S^*) = \text{cov}(\tilde{X}^*|L^*, S^*) + \sigma_X^2 I_{L^*}$, $I_{L^*}$ being the $L^* \times L^*$ identity matrix, that is, the $(j, l)$ entry of the $L^* \times L^*$ matrix $\Sigma_{U^*}$ is $(\Sigma_{U^*})_{j,l} = G_X(S_{ij}, S_{il}) + \sigma_X^2 \delta_{jl}$ with $\delta_{jl} = 1$ if $j = l$ and 0 if $j \neq l$.

According to (10), estimates for the functional principal component scores $\zeta_m^*$ are obtained by substituting estimates of $\mu_X^*$, $\rho_m$ and $\psi_m^*$ that are based on the entire data collection, leading to

$$\hat{\zeta}_m^* = \hat{\rho}_m \hat{\psi}_m^{*T} \widehat{\Sigma}_{U^*}^{-1} (\tilde{U}^* - \hat{\mu}_X^*), \tag{11}$$

where $(\widehat{\Sigma}_{U^*})_{jl} = \widehat{G}_X(S_{ij}, S_{il}) + \hat{\sigma}_X^2 \delta_{jl}$. The predicted trajectory is then obtained as

$$\widehat{Y}_{KM}^*(t) = \hat{\mu}_Y(t) + \sum_{k=1}^{K} \sum_{m=1}^{M} \frac{\hat{\sigma}_{km}}{\hat{\rho}_m} \hat{\zeta}_m^* \hat{\phi}_k(t). \tag{12}$$

In the sparse situation, the Gaussian assumption is crucial. It allows us to obtain the best linear predictors in (10)–(12) through conditional expectations, borrowing strength from the entire sample and thus compensating for the sparseness of data for individual trajectories. Simulation results, reported in Section 4, indicate that the proposed method is quite robust regarding the Gaussian assumption. Theoretical results for predicted trajectories (12) are given in Theorem 2.

3.2. *Asymptotic pointwise confidence bands for response trajectories.* We construct asymptotic confidence bands for the response trajectory $Y^*(\cdot)$ of a new subject, conditional on the sparse and noisy measurements that are available for the underlying predictor function. For $M \geq 1$, let $\zeta^{*M} = (\zeta_1^*, \ldots, \zeta_M^*)^T$, $\tilde{\zeta}^{*M} = (\tilde{\zeta}_1^*, \ldots, \tilde{\zeta}_M^*)^T$, where $\tilde{\zeta}_m^*$ is as in (10), and define the $M \times L^*$ matrix $H = \text{cov}(\zeta^{*M}, \tilde{U}^*|L^*, S^*) = (\rho_1 \psi_1^*, \ldots, \rho_M \psi_M^*)^T$. The covariance matrix of $\tilde{\zeta}^{*M}$ is $\text{cov}(\tilde{\zeta}^{*M}|L^*, S^*) = H \Sigma_{U^*}^{-1} H^T$. Because $\tilde{\zeta}^{*M} = E[\zeta^{*M}|\tilde{U}^*, L^*, S^*]$ is the projection of $\zeta^{*M}$ on the space spanned by the linear functions of $\tilde{U}^*$ given $L^*$ and $S^*$, $\text{cov}(\tilde{\zeta}^{*M} - \zeta^{*M}|L^*, S^*) = \text{cov}(\zeta^{*M}|L^*, S^*) - \text{cov}(\tilde{\zeta}^{*M}|L^*, S^*) \equiv \Omega_M$, where $\Omega_M = D - H \Sigma_{U^*}^{-1} H^T$, and $D = \text{diag}\{\rho_1, \ldots, \rho_M\}$. Under the Gaussian assumption and conditioning on $L^*$ and $S^*$, then $\tilde{\zeta}^{*M} - \zeta^{*M} \sim \mathcal{N}(0, \Omega_M)$.

To construct pointwise bands for $E[Y^*(t)|X^*] = \mu_Y(t) + \sum_{m=1}^{\infty} \sum_{k=1}^{\infty} \sigma_{km} \times \phi_k(t) \tilde{\zeta}_m^*/\rho_m$, let $\widehat{\Omega}_M = \widehat{D} - \widehat{H} \widehat{\Sigma}_{U^*}^{-1} \widehat{H}^T$, where $\widehat{D} = \text{diag}\{\hat{\rho}_1, \ldots, \hat{\rho}_M\}$ and $\widehat{H} = (\hat{\rho}_1 \hat{\psi}_1^*, \ldots, \hat{\rho}_M \hat{\psi}_M^*)^T$. Define $\phi_{tM} = (\phi_1(t), \ldots, \phi_M(t))^T$ for $t \in \mathcal{T}$, and a $K \times M$ matrix $P_{K,M} = (\sigma_{km}/\rho_m)_{1 \leq k \leq K, 1 \leq m \leq M}$. Let $\hat{\phi}_{tK}$, $\widehat{P}_{K,M}$ be the estimates

of $\phi_{tK}$, $P_{KM}$ obtained from the data. Write the prediction in the vector form $\hat{Y}^*_{KM}(t) = \hat{\mu}_Y(t) + \hat{\phi}^T_{tK}\hat{P}_{M,K}\hat{\zeta}^{*M}$. Theorem 3 in Section 6 establishes that the asymptotic distribution of $\{\hat{Y}^*_{KM}(t) - E[Y^*(t)|X^*]\}$ conditional on $L^*$ and $S^*$ can be approximated by $\mathcal{N}(0, \hat{\phi}^T_{tK}\hat{P}_{KM}\hat{\Omega}_M\hat{P}^T_{KM}\hat{\phi}_{tK})$. As a consequence, the $(1-\alpha)$ asymptotic pointwise interval for $E[Y^*(t)|X^*]$, the mean response at predictor level $X^*$, is given by

$$(13) \qquad \hat{Y}^*_{K,M}(t) \pm \Phi(1-\alpha/2)\sqrt{\hat{\phi}^T_{tK}\hat{P}_{KM}\hat{\Omega}_M\hat{P}^T_{KM}\hat{\phi}_{tK}},$$

where $\Phi$ is the standard Gaussian c.d.f.

3.3. *Coefficients of determination for functional linear regression.* In standard linear regression, a measure to quantify the "degree of linear association" between predictor and response variables is the coefficient of determination $R^2$ (e.g., [12], page 138). The coefficient of determination plays an important role in applications of regression analysis, as it may be interpreted as the fraction of variation of the response that is explained by the regression relation. Accordingly, $R^2$ is commonly used as a measure of the strength of the regression relationship.

The proposed extension to functional linear regression can be motivated by the standard linear regression model with a response $Y$ and a predictor $X$, where the coefficient of determination $R^2$ is defined by $R^2 = \text{var}(E[Y|X])/\text{var}(Y)$. This corresponds to the fraction of $\text{var}(Y) = \text{var}(E[Y|X]) + E(\text{var}([Y|X]))$ that is explained by the regression. In the functional setting, the regression function is given by (3), $E[Y(t)|X] = \int_{\mathcal{S}}\beta(s,t)X(s)\,ds$. To measure the global linear association between the functional predictor $X$ and the functional response $Y$, we consider the total variation of $Y$ explained by the regression function to be $\int_{\mathcal{T}}\text{var}(E[Y(t)|X])\,dt$, and by observing (9) and the orthonormality properties of $\{\phi_k\}$ and $\{\psi_m\}$, one obtains

$$(14) \qquad \int_{\mathcal{T}}\text{var}(E[Y(t)|X])\,dt = \sum_{k,m=1}^{\infty}\sigma^2_{km}/\rho_m.$$

The analogous notion of total variation of $Y$ is $\int_{\mathcal{T}}\text{var}(Y(t))\,dt = \int_{\mathcal{T}}G_Y(t,t)\,dt = \sum_{k=1}^{\infty}\lambda_k$. This motivates a functional version of $R^2$,

$$(15) \qquad R^2 = \frac{\int_{\mathcal{T}}\text{var}(E[Y(t)|X])\,dt}{\int_{\mathcal{T}}\text{var}(Y(t))\,dt} = \frac{\sum_{k,m=1}^{\infty}\sigma^2_{km}/\rho_m}{\sum_{k=1}^{\infty}\lambda_k}.$$

Since $\text{var}(E[Y|X]) \leq \text{var}(Y)$ for random variables $X, Y$, it follows that $\int_{\mathcal{T}}\text{var}(E[Y(t)|X])\,dt \leq \int_{\mathcal{T}}\text{var}(Y(t))\,dt$, that is, $\sum_{m,k=1}^{\infty}\sigma^2_{mk}/\rho_m \leq \sum_{k=1}^{\infty}\lambda_k < \infty$. Thus, the functional $R^2$ (15) always satisfies $0 \leq R^2 \leq 1$.

Another interpretation of the functional coefficient of determination $R^2$ (15) is as follows: Denoting by $R^2_{km}$ the coefficients of determination of the simple linear regressions of the functional principal component scores $\xi_k$ on $\zeta_m$, $1 \leq k, m < \infty$, one finds $R_{km}^2 = \sigma_{km}^2/(\rho_m \lambda_k)$, and $R_k^2 = \sum_{m=1}^{\infty} R_{km}^2$ is the coefficient of determination of regressing $\xi_k$ on all $\zeta_m$, $m = 1, 2, \ldots$, simultaneously, by successively adding predictors $\zeta_m$ into the regression equation and observing that $R_k^2$ is obtained as the sum of the $R_{km}^2$, as the predictors $\zeta_m$ are uncorrelated. Then $R^2 = \sum_{k=1}^{\infty} \lambda_k R_k^2/(\sum_{k=1}^{\infty} \lambda_k)$ is seen to be a weighted average of these $R_k^2$, with weights provided by the $\lambda_k$. According to (15), a natural estimate $\widehat{R}^2$ for the functional coefficient of determination $R^2$ is

$$(16) \qquad \widehat{R}^2 = \frac{\sum_{k=1}^{K} \sum_{m=1}^{M} \hat{\sigma}_{km}^2/\hat{\rho}_m}{\sum_{k=1}^{K} \hat{\lambda}_k},$$

where $\hat{\sigma}_{km}$ are as in (7).

Besides the functional $R^2$ (15) that provides a global measure of the linear association between processes $X$ and $Y$, we also propose a local version of a functional coefficient of determination. The corresponding function $R^2(t)$ may be considered a functional extension of the local $R^2$ measure that has been introduced in [10] and [11]. As shown above, for fixed $t \in \mathcal{T}$, the variation of $Y(t)$ explained by the predictor process $X$ is determined by $\text{var}(E[Y(t)|X])/\text{var}(Y(t))$. This motivates the following definition of a pointwise functional coefficient of determination $R^2(t)$:

$$(17) \quad R^2(t) = \frac{\text{var}(E[Y(t)|X])}{\text{var}(Y(t))} = \frac{\sum_{m=1}^{\infty} \sum_{k,\ell=1}^{\infty} \sigma_{km} \sigma_{\ell m} \phi_k(t) \phi_\ell(t)/\rho_m}{\sum_{k=1}^{\infty} \lambda_k \phi_k^2(t)}.$$

Note that $R^2(t)$ satisfies $0 \leq R^2(t) \leq 1$ for all $t \in \mathcal{T}$.

A second option to obtain an overall $R^2$ value is to extend the pointwise measure $R^2(t)$ to a global measure by taking its integral, which leads to an alternative definition of the global functional coefficient of determination, denoted by $\widetilde{R}^2$,

$$
\begin{aligned}
(18) \qquad \widetilde{R}^2 &= \frac{1}{|\mathcal{T}|} \int_{\mathcal{T}} \frac{\text{var}(E[Y(t)|X])}{\text{var}(Y(t))} \, dt \\
&= \frac{1}{|\mathcal{T}|} \int_{\mathcal{T}} \frac{\sum_{m=1}^{\infty} \sum_{k,\ell=1}^{\infty} \sigma_{km} \sigma_{\ell m} \phi_k(t) \phi_\ell(t)/\rho_m}{\sum_{k=1}^{\infty} \lambda_k \phi_k^2(t)} \, dt,
\end{aligned}
$$

where $|\mathcal{T}|$ denotes the length of the time domain $\mathcal{T}$. Natural estimates of $R^2(t)$ and $\widetilde{R}^2$ are given by

$$(19) \qquad \widehat{R}^2(t) = \frac{\sum_{m=1}^{M} \sum_{k,\ell=1}^{K} \hat{\sigma}_{km} \hat{\sigma}_{\ell m} \hat{\phi}_k(t) \hat{\phi}_\ell(t)/\hat{\rho}_m}{\sum_{k=1}^{K} \hat{\lambda}_k \hat{\phi}_k^2(t)},$$

$$(20) \qquad \widehat{\widetilde{R}}^2 = \frac{1}{|\mathcal{T}|} \int_{\mathcal{T}} \frac{\sum_{m=1}^{M} \sum_{k,\ell=1}^{K} \hat{\sigma}_{km} \hat{\sigma}_{\ell m} \hat{\phi}_k(t) \hat{\phi}_\ell(t)/\hat{\rho}_m}{\sum_{k=1}^{K} \hat{\lambda}_k \hat{\phi}_k^2(t)} \, dt.$$

We refer to Section 5 for further discussion of $R^2$, $R^2(t)$ and $\widetilde{R}^2$ in applications and to Theorem 4 in Section 6 regarding the asymptotic convergence of these estimates.

**4. Simulation studies.** Simulation studies were based on 500 i.i.d. normal and 500 i.i.d. mixture samples, where 100 pairs of response and predictor trajectories were generated for each sample. Emulating very sparse and irregular designs, each pair of response and predictor functions was observed at different sets of time points. The number of measurements was randomly chosen for each predictor and each response trajectory, with equal probability from $\{3, 4, 5\}$ for $X_i$ uniformly, and independently chosen from $\{3, 4, 5\}$ for $Y_i$ uniformly, also with equal probability. This setup reflects very sparse designs with at most five observations available per subject. Once their numbers were determined, the locations of the measurements were uniformly distributed on $[0, 10]$ for both $X_i$ and $Y_i$, respectively.

The predictor trajectories $X_i$ and associated sparse and noisy measurements $U_{il}$ (1) were generated as follows. The simulated processes $X$ had the mean function $\mu_X(s) = s + \sin(s)$, with covariance function constructed from two eigenfunctions, $\psi_1(s) = -\cos(\pi s/10)/\sqrt{5}$ and $\psi_2(s) = \sin(\pi s/10)/\sqrt{5}$, $0 \le s \le 10$. We chose $\rho_1 = 2$, $\rho_2 = 1$ and $\rho_m = 0$, $m \ge 3$, as eigenvalues, and $\sigma_X^2 = 0.25$ as the variance of the additional measurement errors $\varepsilon_{il}$ in (1), which were assumed to be normal with mean 0. For the 500 normal samples, the FPC scores $\zeta_{im}$ ($m = 1, 2$) were generated from $\mathcal{N}(0, \rho_m)$, while the $\zeta_{im}$ for the nonnormal samples were generated from a mixture of two normals, $\mathcal{N}(\sqrt{\rho_m/2}, \rho_m/2)$ with probability $1/2$ and $\mathcal{N}(-\sqrt{\rho_m/2}, \rho_m/2)$, also with probability $1/2$.

For the response trajectories, letting $b_{11} = 2$, $b_{12} = 2$, $b_{21} = 1$, $b_{22} = 2$, the regression function was $\beta(s, t) = \sum_{k=1}^{2} \sum_{m=1}^{2} b_{km} \psi_m(s) \psi_k(t)$, $t, s \in [0, 10]$, and the response trajectories were $E[Y_i(t)|X_i] = \int_0^{10} \beta(s, t) X_i(s)\, ds$. Only the sparse and noisy observations $V_{ij} = E[Y_i(T_{ij})|X_i] + \epsilon_{ij}$ (2) were available for response trajectories, contaminated with pseudo-i.i.d. errors $\epsilon_{ij}$ with density $\mathcal{N}(0, 0.1)$.

We investigated predicting response curves for new subjects. For each Monte Carlo simulation run, we generated 100 new predictor curves $X_i^*$, with noisy measurements taken at the same random time points as $X_i$, and 100 associated response curves $E[Y_i^*|X_i^*]$. Relative mean squared prediction error was used as an evaluation criterion, given by

$$(21) \qquad \text{RMSPE} = \frac{1}{n} \sum_{i=1}^{n} \int_0^{10} \frac{\{\widehat{Y}_{i,KM}^*(t) - E[Y_i^*(t)|X_i^*]\}^2\, dt}{\int_0^{10} \{E[Y_i^*(t)|X_i^*]\}^2\, dt},$$

where predicted trajectories $\widehat{Y}_{i,KM}^*$ were obtained according to (11) and (12). This method is denoted as CE in Table 1, and was compared with a "classical" functional regression approach that was also based on (12), but with the conditional expectation replaced by the integral approximation $\hat{\zeta}_{im}^{*I} = \sum_{l=1}^{L_i^*}(U_{il}^* - \hat{\mu}_X(S_{il}^*))\hat{\psi}_m(S_{il}^*)(S_{il}^* - S_{i,l-1}^*)$, denoted by IN in Table 1. The numbers of eigenfunctions $K$ and $M$ (35) were chosen by the AIC criterion (35), separately for each simulation. We also included the case of irregular but nonsparse data, where the random numbers of the repeated measurements were chosen from $\{20, \ldots, 30\}$ for both $X_i$ and $Y_i$ with equal probability. From the results in Table 1, we see that, for sparse data, the CE method improves the prediction errors by 57%/60% for normal/mixture distributions, while the gains for nonsparse data are not as dramatic, but nevertheless present. The CE method emerges as superior for the sparse data case.

## 5. Applications.

5.1. *Primary biliary cirrhosis data.* Primary biliary cirrhosis [23] is a rare but fatal chronic liver disease of unknown cause, with a prevalence of about 50 cases per million population. The data were collected between January 1974 and May 1984 by the Mayo Clinic (see also Appendix D of [17]). The patients were scheduled to have measurements of blood characteristics at six months, one year and annually thereafter post diagnosis. However, since many individuals missed some of their scheduled visits, the data are sparse and irregular with unequal numbers of repeated measurements per subject and also different measurement times $T_{ij}$ per individual.

To demonstrate the usefulness of the proposed methods, we explore the dynamic relationship between albumin in mg/dl (predictor) and prothrombin time in seconds (response), which are both longitudinally measured. We include 137 female patients, and the measurements of albumin level and prothrombin time before 2500 days. For both albumin and prothrombin time, the number of observations ranged from 1 to 10, with a median of 5 measurements per subject. Individual trajectories of albumin and prothrombin time are shown in Figure 1.

TABLE 1
*Results of 500 Monte Carlo runs with $n = 100$ trajectories per sample. Shown are medians of observed squared prediction errors, RMSPE (21). Here CE is the conditional expectation method (11), (12) and IN stands for integral approximation*

| Regression | | Normal | Mixture |
|---|---|---|---|
| Sparse | CE | 0.0083 | 0.0081 |
| | IN | 0.0193 | 0.0204 |
| Nonsparse | CE | 0.0057 | 0.0062 |
| | IN | 0.0078 | 0.0079 |

The smooth estimates of the mean function for both albumin and prothrombin time are also displayed in Figure 1, indicating opposite trends. The AIC criterion leads to the choice of two eigenfunctions for both predictor and response functions, and smooth eigenfunction estimates are presented in Figure 2. For both albumin and prothrombin time, the first eigenfunction reflects an overall level, and the second eigenfunction a contrast between early and late times. The estimate of the regression function $\beta$ is displayed in Figure 3. Its shape implies that, for the prediction of early prothrombin times, late albumin levels contribute positively, while early levels contribute negatively, whereas the prediction of late prothrombin times is based on a sharper contrast with highly positive weighting of early albumin levels and negative weighting of later levels.

We randomly selected four patients from the sample for which the trajectory of prothrombin time was to be predicted solely from the sparse and noisy albumin measurements. For this prediction, the data of each subject to be predicted were omitted, the functional regression model was fitted from the remaining subjects, and then the predictor measurements were entered into the fitted model to obtain the predicted response trajectory, thus leading to genuine predictions. Predicted curves and 95% pointwise confidence bands are shown in Figure 4. Note that these predictions of longitudinal trajectories are based on just a few albumin measurements, and the prothrombin time response measurements shown in the figures are not used.

Regarding the functional coefficients of determination $R^2$ (15) and $\widetilde{R}^2$ (18), we obtain very similar estimates, $\widehat{R}^2 = 0.37$ and $\widehat{\widetilde{R}}^2 = 0.36$, which we would interpret to mean that about 36% of the total functional variation of the prothrombin time trajectories is explained by the albumin data, indicating a reasonably strong functional regression relationship. The curve of estimated pointwise coefficients of determination $\widehat{R}^2(t)$ (19) is shown in Figure 5, left panel, which describes the trend of the proportion of the variation of the prothrombin time, at each argument value, that is explained by the entire albumin trajectories. We find that the observations in the second half are better determined by the albumin trajectories than the values in the first half of the domain of prothrombin time.

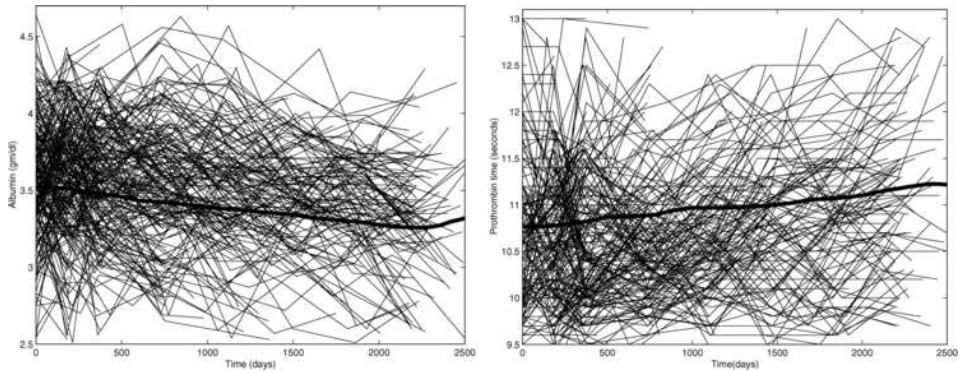

Fig. 1.  *Left panel: Observed individual trajectories (solid) and the smooth estimate of the mean function for albumin (thick solid). Right panel: Corresponding observed individual trajectories (solid) and the smooth estimate of the mean function for prothrombin time (thick solid), for the primary biliary cirrhosis (PBC) data.*

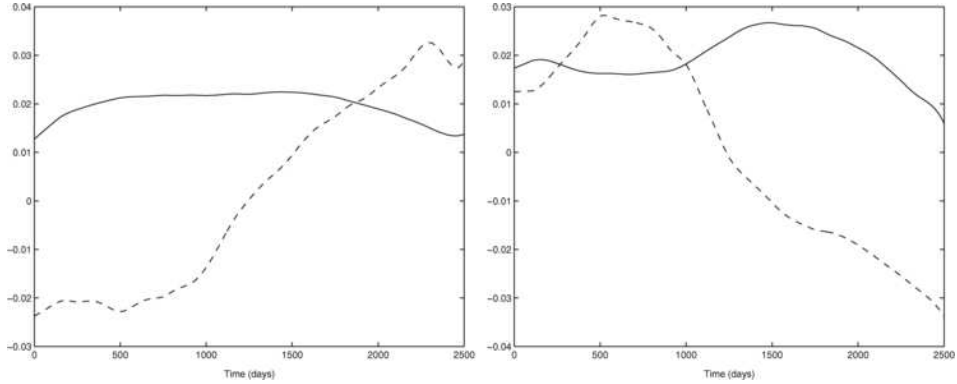

FIG. 2. *Left panel: Smooth estimates of the first (solid) and second (dashed) eigenfunctions for albumin, accounting for 87% and 8% of total variation. Right panel: Smooth estimates of the first (solid) and second (dashed) eigenfunctions for prothrombin time, accounting for 54% and 33% of total variation, for the PBC data.*

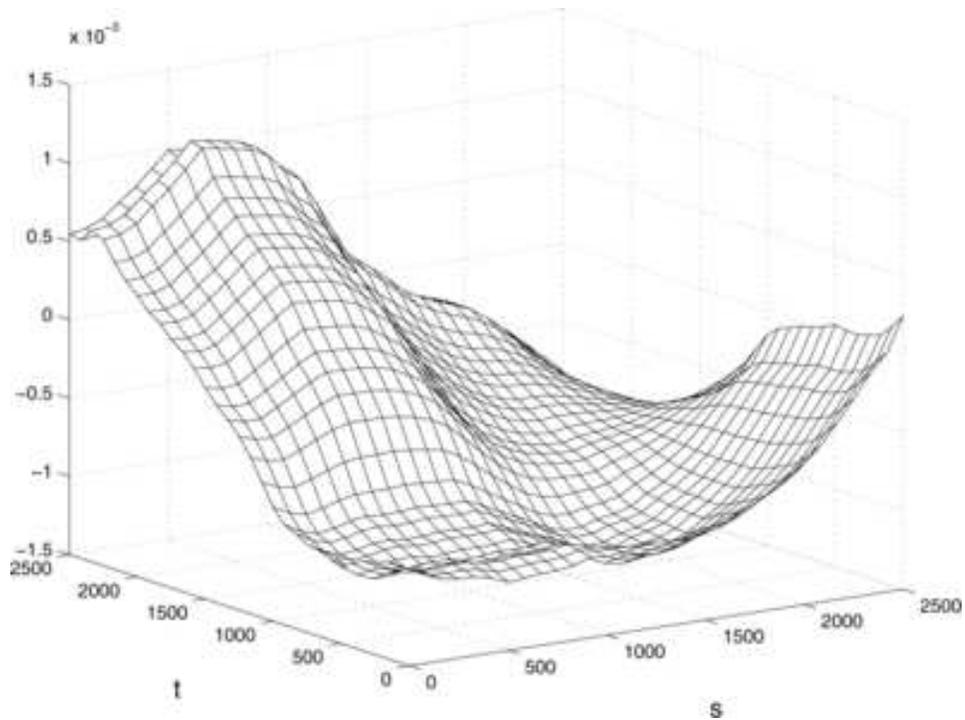

FIG. 3. *Estimated regression function* (8), *where the predictor (albumin) time is s (in days), and the response (prothrombin) time is t (in days), for the PBC data.*

5.2. *Functional regression of systolic blood pressure on body mass index.* As a second example, we discuss a functional regression analysis of systolic blood pressure trajectories (responses) on body mass index trajectories (predictor), using anonymous data from the Baltimore Longitudinal Study of Aging (BLSA), a major longitudinal study of human aging [29]. The data consist of 1590 male volunteers who were scheduled to visit the Gerontology

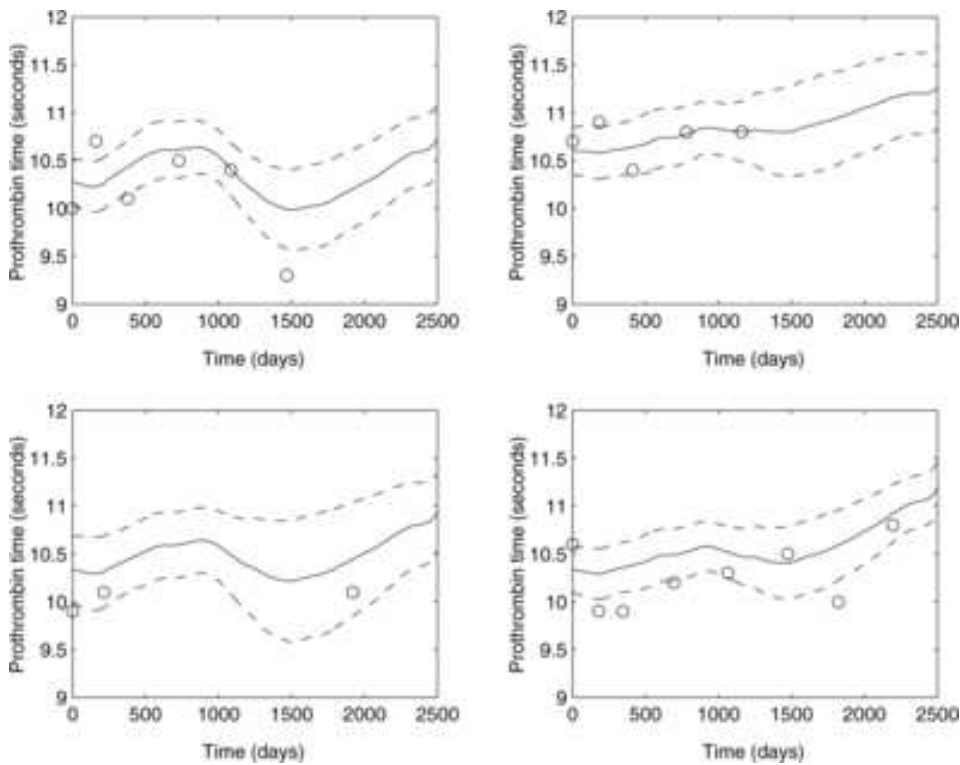

Fig. 4. *Observed values (circles) for prothrombin times (not used for prediction), predicted curves (solid) and 95% pointwise bands (dashed), for four randomly selected patients, where bands and predicted curves are based on one-curve-leave-out analysis and sparse noisy albumin measurements (not shown) of the predictor function.*

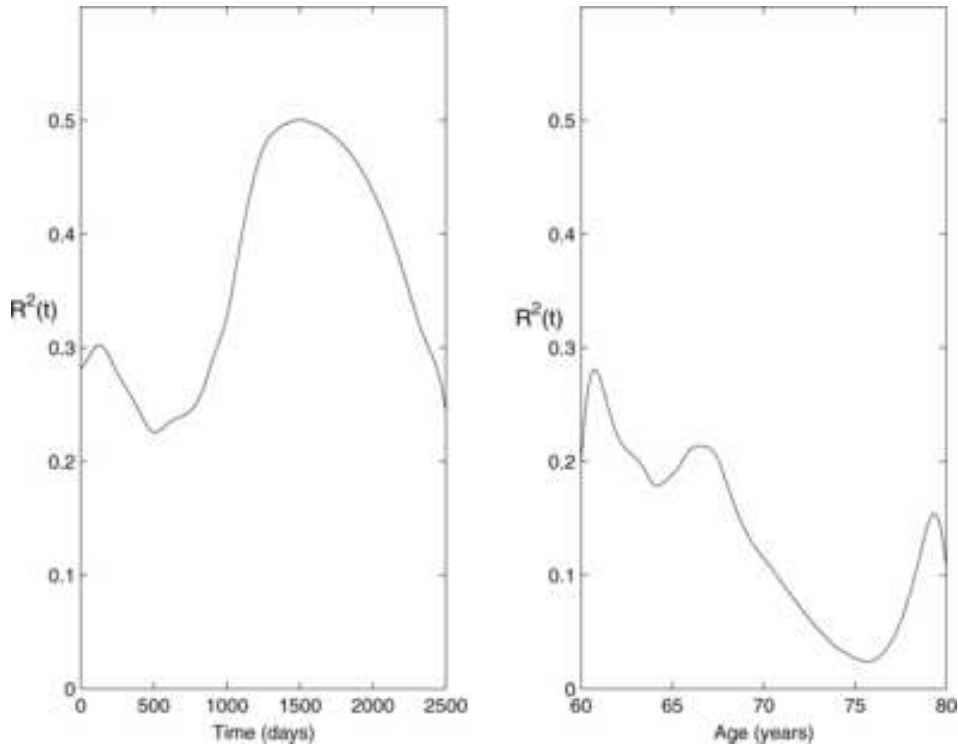

Fig. 5. *Estimated pointwise coefficient of determination $\widehat{R}^2(t)$ (19) as a function of time t for the PBC data (left panel, time is duration of study in days) and for the BLSA data (right panel, time is age in years).*

Research Center bi-annually. Time corresponds to age of each subject and is measured in years. On each visit, systolic blood pressure (in mm Hg) and body mass index (BMI = weight in kg/height in $m^2$) were assessed. Since both measurements are highly variable, the data are noisy, and as many participants missed scheduled visits, or were seen at other than the scheduled times, the data are sparse with largely unequal numbers of repeated measurements and widely differing measurement times per subject. More details about the study and data can be found in [24], and for previous statistical approaches, we refer to [22].

We included the participants with measurements within the age range [60, 80], and first checked for outliers based on standardized residuals of body mass index (BMI) and systolic blood pressure (SBP), respectively. The standardized residuals are defined as residuals divided by the pooled sample standard deviation, where residuals are the differences between measurements and the estimated mean function obtained by scatterplot smoothing, using the local linear smoother. We excluded subjects with standardized residuals larger (or less) than ±3, for either BMI or SBP. Individual trajectories of BMI and SBP for the included 812 subjects are shown in Figure 6, along with the smooth estimated mean functions of BMI and SBP.

While average BMI decreases after age 64, SBP throughout shows an increasing trend. Based on the AIC criterion, three eigenfunctions are used for the predictor (BMI) function, and four for the response (SBP) function; these are displayed in Figure 7. The first eigenfunctions of both processes correspond to an overall mean effect, and the second eigenfunctions to a contrast between early and late ages, with further oscillations reflected in third and fourth eigenfunctions.

The estimated regression function in Figure 8 indicates that a contrast between late and early BMI levels forms the prediction of SBP levels at later ages, where late BMI levels are weighted positive and early levels negative. When predicting SBP at age 80, the entire BMI trajectory matters, and rapid overall declines in BMI lead to the lowest SBPs, where speed of de-cline between 60 and 65 and between 75 and 80 is critical. Similar patterns can be identified for predicting SBP at other ages. As in the previous example, we randomly select four study participants and obtain predictions and 95% pointwise bands for each of these, based on one-leave-out functional regression analysis (Figure 9). The predicted trajectories are found to be reasonably close to the observations, which are not used in the analysis.

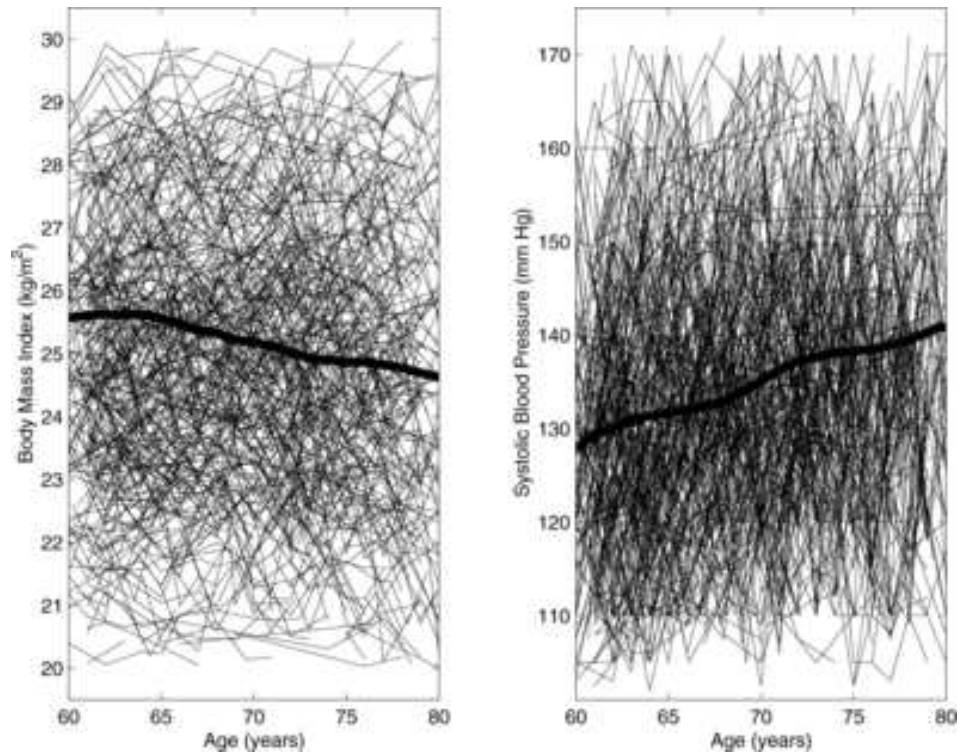

FIG. 6. *Left panel: Observed individual trajectories (solid) and the smooth estimate of the mean function (thick solid) for body mass index (BMI). Right panel: Corresponding observed individual trajectories (solid) and the smooth estimate of the mean function (thick solid) for systolic blood pressure (SBP), for the BLSA data.*

The functional coefficients of determination $R^2$ (15) and $\widetilde{R}^2$ (18) were both estimated as 0.13, indicating that the dynamics of body mass index explains 13% of the total variation of systolic blood pressure trajectories; the functional regression relationship is seen to be weaker than in the previous example. The curve of estimated pointwise coefficients of determination $R^2(t)$ (17) is displayed in Figure 5, right panel, indicating generally weaker linear association at older ages (beyond 70 years) as compared to the earlier ages (60 to 70 years). The minimal linear association between predictor trajectories and the functional response is seen to occur around age 75.7.

**6. Asymptotic properties.** In this section we present the consistency of the regression function estimate (8) in Theorem 1. The consistency of predicted trajectories is reflected in Theorem 2, and the asymptotic distributions of predicted trajectories in Theorem 3. The proposed functional coefficient of determination $\widehat{R}^2$ (16) is shown to be consistent for $R^2$ (15) in Theorem 4.

In what follows, we only consider the case that the processes $X$ and $Y$ are infinite-dimensional. If they are finite-dimensional and there exist true finite $K$ and $M$, such that $G_X$ and $G_Y$ are finite-dimensional surfaces, then, as $n$ tends to infinity, remainder terms such as $\theta_n$ and $\vartheta_n$ in (41) will disappear. Appropriately modified versions of the following theoretical results still hold. Since in this case the true $K$ and $M$ are finite, most of the technical assumptions that are needed to handle the infinite case would not be needed, such as (A1)–(A3), (A6)–(A7) and (B5), and $K(n), M(n)$ would be assumed to converge to the true $K$ and $M$ instead of infinity, as $n \to \infty$.

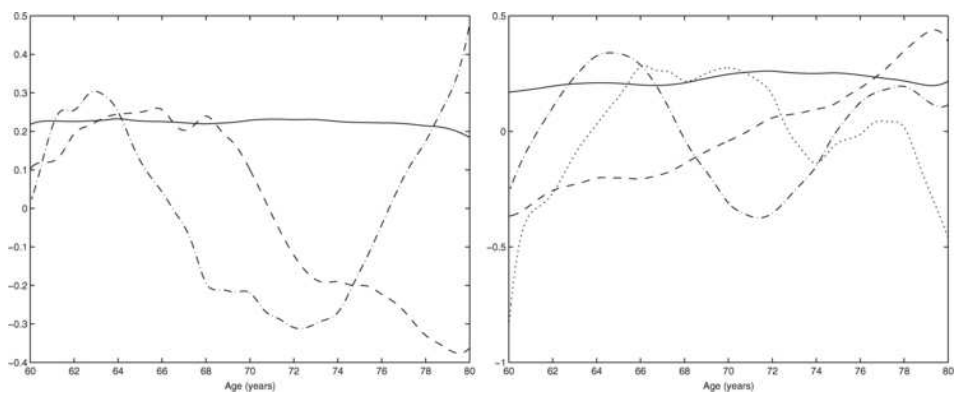

FIG. 7. *Left panel: Smooth estimates of the first (solid), second (dashed) and third (dash-dot) eigenfunctions for BMI, accounting for 90%, 6% and 3% of total variation. Right panel: Smooth estimates of the first (solid), second (dashed), third (dash-dot) and fourth (dotted) eigenfunctions for SPB, accounting for 76%, 15%, 4% and 2% of total variation, for the BLSA data.*

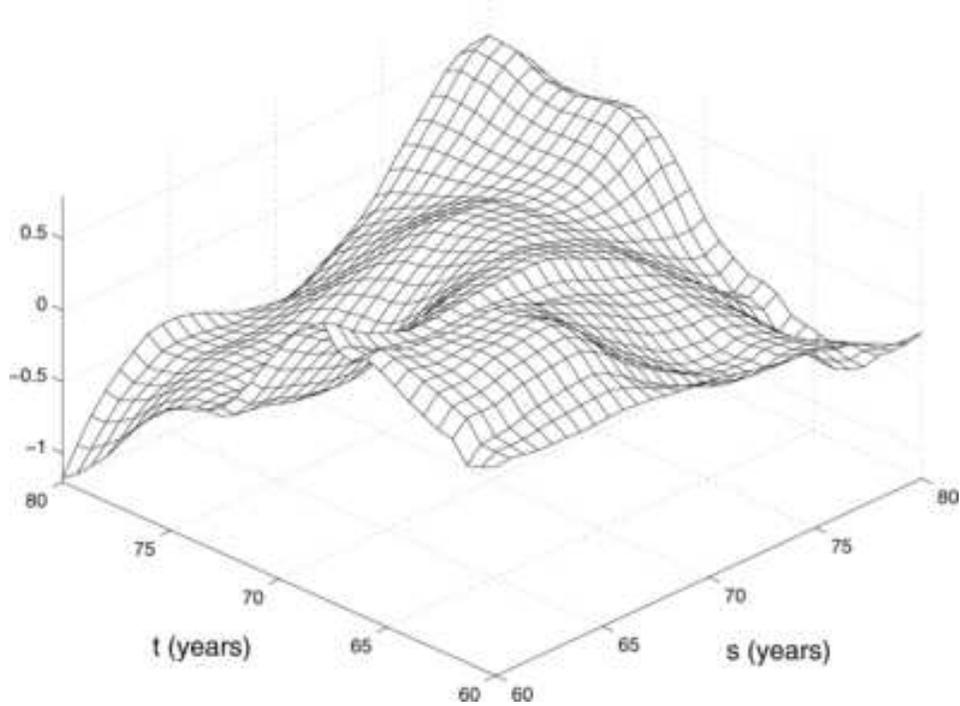

Fig. 8.   *Estimated regression function* (8), *where the predictor (BMI) time is s (in years), and the response (SBP) time is t (in years), for the BLSA data.*

To define the convergence of the right-hand side of (5) in the $L^2$ sense, we require that

(A1)  $\sum_{k=1}^{\infty} \sum_{m=1}^{\infty} \sigma_{km}^2 / \rho_m^2 < \infty$.

Furthermore, the right-hand side of (5) converges uniformly on $\mathcal{S} \times \mathcal{T}$, provided that

(A2)  $\gamma(s,t) = \sum_{k=1}^{\infty} \sum_{m=1}^{\infty} |\sigma_{km} \psi_m(s) \phi_k(t)| / \rho_m$ is continuous in $s$ and $t$, and the function $\beta_{KM}(s,t) = \sum_{k=1}^{K} \sum_{m=1}^{M} \sigma_{km} \psi_m(s) \phi_k(t) / \rho_m$ absolutely converges to $\beta(s,t)$ for all $s \in \mathcal{S}$, $t \in \mathcal{T}$ as $M, K \to \infty$.

The numbers $M = M(n)$ and $K = K(n)$ of included eigenfunctions are integer-valued sequences that depend on the sample size $n$; see assumption (B5) in Appendix A.1. For simplicity, we suppress the dependency of $M$ and $K$ on $n$ in the notation. The consistency of $\hat{\beta}$ (8) is obtained as follows.

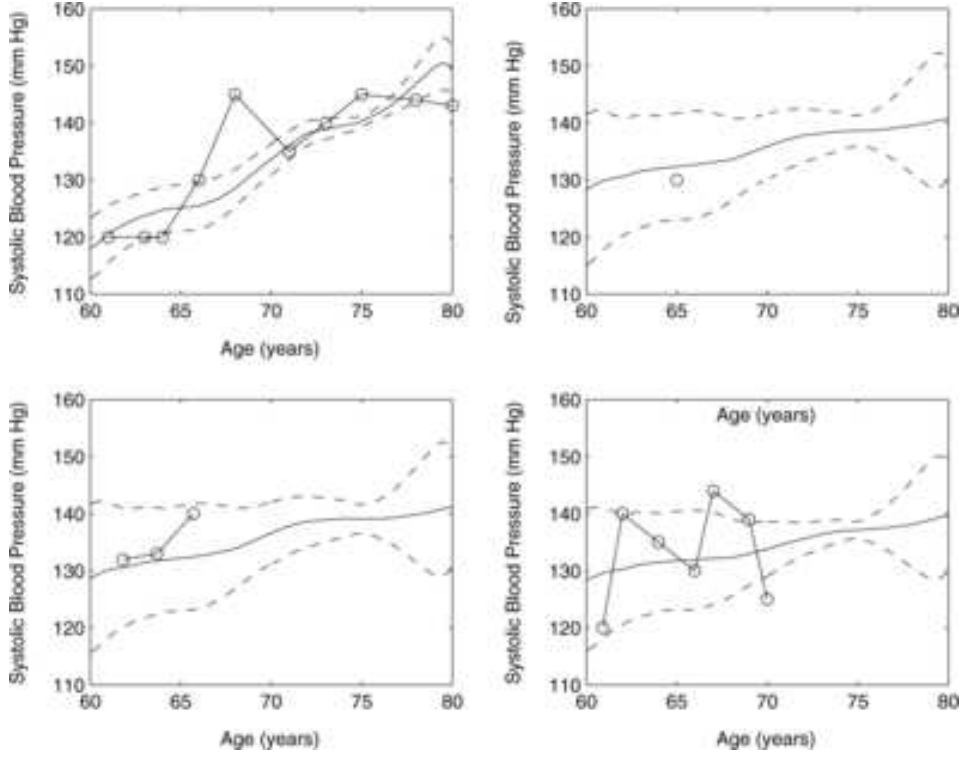

Fig. 9. *Observed values (circles) of systolic blood pressure (SBP) that are not used for the prediction, predicted curves (solid) and 95% pointwise bands (dashed), for four randomly selected patients, where bands and predicted curves are based solely on sparse and noisy measurements of the predictor function (BMI, not shown).*

Theorem 1. *Under* (A1) *and the assumptions of Lemma* A.1 *and* (B5) *(see Appendix* A.1*),*

$$\lim_{n\to\infty} \int_{\mathcal{T}} \int_{\mathcal{S}} [\hat{\beta}(s,t) - \beta(s,t)]^2 \, ds \, dt = 0 \qquad in \ probability, \tag{22}$$

*and if* (A1) *is replaced with* (A2),

$$\lim_{n\to\infty} \sup_{(s,t)\in\mathcal{S}\times\mathcal{T}} |\hat{\beta}(s,t) - \beta(s,t)| = 0 \qquad in \ probability. \tag{23}$$

The rate of convergence in (22) and (23) depends on specific properties of processes $X$ and $Y$ in the following way: If $\tau_n, v_n$ and $\varsigma_n$ are defined as in (B.5) of Appendix A.1, and $\theta_n$, $\vartheta_n$ are defined as in (41) in Appendix A.2, then, for (22) we obtain the rate

$$\int_{\mathcal{T}} \int_{\mathcal{S}} [\hat{\beta}(s,t) - \beta(s,t)]^2 \, ds \, dt = O_p(\tau_n + \upsilon_n + \varsigma_n + \theta_n),$$

and for (23) the rate is

$$\sup_{(s,t) \in \mathcal{S} \times \mathcal{T}} |\hat{\beta}(s,t) - \beta(s,t)| = O_p(\tau_n + \upsilon_n + \varsigma_n + \vartheta_n),$$

as $n \to \infty$. Here $\varsigma_n$ depends on bandwidths $h_1$ and $h_2$ that are used in the smoothing step (36) for the cross-covariance function $C(s,t) = \operatorname{cov}(X(s), Y(t))$. These rates depend on specific properties of the processes $X$ and $Y$, such as the spacings of the eigenvalues of their autocovariance operators. We note that, due to the sparsity of the data (at most, finitely many observations are made per random trajectory), fast rates of convergence cannot be expected in this situation, in contrast to the case where entire trajectories are observed or are densely sampled.

Recall that $Y_{KM}^*(t) = \mu_Y(t) + \sum_{k=1}^{K} \sum_{m=1}^{M} \sigma_{km} \phi_k(t) \zeta_m^* / \rho_m$, and the prediction $\hat{Y}_{KM}^*(t) = \hat{\mu}_Y(t) + \sum_{k=1}^{K} \sum_{m=1}^{M} \hat{\sigma}_{km} \hat{\phi}_k(t) \hat{\zeta}_m^* / \hat{\rho}_m$, where $\hat{\zeta}_m^*$ is as in (11). Define the target trajectory $\tilde{Y}^*(t) = \mu(t) + \sum_{k=1}^{\infty} \sum_{m=1}^{\infty} \sigma_{km} \tilde{\zeta}_m^* \phi_k(t) / \rho_m$ and $\tilde{Y}_{KM} = \mu(t) + \sum_{k=1}^{K} \sum_{m=1}^{M} \sigma_{km} \tilde{\zeta}_m^* \phi_k(t) / \rho_m$, where $\tilde{\zeta}_m^*$ is defined in (10). Assume that

(A3) $\sum_{k=1}^{\infty} \sum_{m=1}^{\infty} \sigma_{km}^2 / (\lambda_k \rho_m) < \infty$.

Then $E[Y^*(t)|X^*]$ and $\tilde{Y}^*$ are defined as the limits of $Y_{KM}^*(t)$ and $\tilde{Y}_{KM}^*(t)$ in the $L^2$ sense; see Lemma A.3 in Appendix A.3. Furthermore, we assume:

(A4) The number and locations of measurements for a given subject or cluster remain unaltered as the sample size $n \to \infty$.

THEOREM 2. *Under* (A3), (A4) *and the assumptions of Lemma* A.1 *and* (B5) *(see Appendix* A.1*), given* $L^*$ *and* $S^*$, *for all* $t \in \mathcal{T}$,

$$(24) \qquad \lim_{n \to \infty} \hat{Y}_{KM}^*(t) = \tilde{Y}^*(t) \qquad \text{in probability.}$$

This provides the consistency of the prediction $\hat{Y}_{KM}^*$ for the target trajectory $\tilde{Y}^*$.

For the following results, we require Gaussian assumptions.

(A5) For all $1 \leq i \leq n$, $m \geq 1$ and $1 \leq l \leq L_i$, the functional principal component scores $\zeta_{im}$ and the measurement errors $\varepsilon_{il}$ in (1) are jointly Gaussian.

Define $\omega_{KM}(t_1, t_2) = \phi_{t_1 K}^T P_{KM} \Omega_M P_{KM}^T \phi_{t_2 K}$, for $t_1, t_2 \in \mathcal{T}$. Then $\omega_{KM}(t_1, t_2)$ is a sequence of continuous positive definite functions, and $\hat{\omega}_{KM}(t_1, t_2) = \hat{\phi}_{t_1 K}^T \hat{P}_{KM} \hat{\Omega}_M \hat{P}_{KM}^T \hat{\phi}_{t_2 K}$ is an estimate of $\omega_{KM}(t_1, t_2)$. We require existence of a limiting function and, therefore, the following analytical condition:

(A6) There exists a continuous positive definite function $\omega(t_1, t_2)$ such that $\omega_{KM}(t_1, t_2) \to \omega(t_1, t_2)$ as $K, M \to \infty$.

We obtain the asymptotic distribution of $\{\hat{Y}_{KM}^*(t) - E[Y^*(t)|X^*]\}$ as follows, providing inference for predicted trajectories.

THEOREM 3. *Under* (A3)–(A6) *and the assumptions of Lemma* A.1 *and* (B5) *(see Appendix* A.1*), given $L^*$ and $S^*$, for all $t \in \mathcal{T}$, $x \in \Re$,*

$$(25) \qquad \lim_{n \to \infty} P\left\{ \frac{\hat{Y}_{KM}^*(t) - E[Y^*(t)|X^*]}{\sqrt{\hat{\omega}_{KM}(t,t)}} \le x \right\} = \Phi(x).$$

Considering the measure $R^2$ (15), $R^2$ is well defined since $\sum_{m,k=1}^{\infty} \sigma_{mk}^2 / \rho_m \le \sum_{k=1}^{\infty} \lambda_k < \infty$. Furthermore, the right-hand side of (17) uniformly converges on $t \in \mathcal{T}$, provided that

(A7) $\kappa(t) = \sum_{m=1}^{\infty} \sum_{k,\ell=1}^{\infty} |\sigma_{km} \sigma_{\ell m} \phi_k(t) \phi_\ell(t)| / \rho_m$ is continuous in $t \in \mathcal{T}$, and the function $R_{MK}^2(t) = \sum_{m=1}^{M} \sum_{k,\ell=1}^{K} \sigma_{km} \sigma_{\ell m} \phi_k(t) \phi_\ell(t) / \rho_m$ absolutely converges to $R^2(t)$ for all $t \in \mathcal{T}$.

The consistency of $\hat{R}^2$ (16), $\hat{R}^2(t)$ (19) and $\hat{\tilde{R}}^2$ (20) is obtained as a consequence of Lemma A.1, (A7) and (B5).

THEOREM 4. *Under the assumptions of Lemma* A.1 *and* (B5) *(see Appendix* A.1*),*

$$(26) \qquad \lim_{n \to \infty} \hat{R}^2 = R^2 \qquad \text{in probability,}$$

*and if* (A7) *is assumed,*

$$(27) \qquad \lim_{n \to \infty} \sup_{t \in \mathcal{T}} |\hat{R}^2(t) - R^2(t)| = 0 \qquad \text{in probability,}$$

$$(28) \qquad \lim_{n \to \infty} \hat{\tilde{R}}^2 = \tilde{R}^2 \qquad \text{in probability.}$$

We note that the rate of convergence in (26) is the same as in the remark after Theorem 1, and that the rate in (27) and (28) is given by $O_p(\tau_n + \upsilon_n + \varsigma_n + \pi_n)$, where $\pi_n$ is as in (29).

**7. Concluding remarks.**   The functional regression method we are proposing applies to the situation where both predictors and responses are curves. Sparse data situations will occur also in other functional regression situations where responses could be functions and predictors vectors, or where responses are scalars and the predictors are curves. The ideas presented here can be extended to such situations.

The approach to functional regression we have proposed is quite flexible, and in simulations is seen to be robust to violations of assumptions such as the Gaussian assumption. It is a useful tool for data where both predictors and responses are contaminated by errors. We refer to [9] for another approach and discussion of de-noising of single pairs of curves observed with measurement errors. Besides varying coefficient models, the available methodology for situations where one has a sample of predictor and response functions is quite limited. Common methods in longitudinal data analysis such as Generalized Estimating Equations and Generalized Linear Mixed Models are not suitable for this task.

The proposed methodology may prove generally useful for longitudinal data with missing measurements, where missingness would be assumed to be totally unrelated to the random trajectories and errors. Extensions to situations where missingness is correlated with the time courses would be of interest in many practical applications. There are also some limitations to the functional regression approach under sparse data. Our prediction methods target the trajectory conditional on the available data, while the response trajectory given the entire but unobservable predictor trajectory is not accessible. While in theory it is enough that the probability that one observes more than one measurement per random trajectory is positive, in practice there needs to be a substantial number of subjects with two or more observations.

Sometimes a prediction for a response trajectory may be desired even if there is no observation at all available for that subject. In this case we predict the estimated mean response function as the response trajectory. This is not unreasonable, since borrowing strength from other subjects to predict the response trajectory for a given subject is a key feature of the proposed method that will come more into play for subjects with very few measurements. Predictions for these subjects will often be relatively closer to the mean response than for subjects with many measurements.

We conclude by remarking that extensions to cases with more than one predictor function are of interest in a number of applications, and would be analogous to the extension of simple linear regression to multiple linear regression. Functional regression is only at its initial stages and much more work needs to be done.

## APPENDIX

**A.1. Assumptions and notation.** The data $(S_{il}, U_{il})$ and $(T_{ij}, V_{ij})$, $i = 1, \ldots, n$, $l = 1, \ldots, L_i$, $j = 1, \ldots, N_i$, as described in (1) and (2), are assumed to have the same distribution as $(S, U)$ and $(T, V)$, with joint densities $g_1(s, u)$ and $g_2(t, v)$. Assume that the observation times $S_{il}$ are i.i.d. with marginal density $f_S(s)$. Dependence is allowed for observations $U_{il_1}$ and $U_{il_2}$ made for the same subject or cluster, and analogous properties hold for $V_{ij}$, where $T_{ij}$ are i.i.d. with marginal density $f_T(t)$. We make the following assumptions for the number of observations $L_i$ and $N_i$ that are available for the $i$th subject:

(B1.1) The number of observations $L_i$ and $N_i$ made for the $i$th subject or cluster are random variables such that $L_i \overset{\text{i.i.d.}}{\sim} L$, $N_i \overset{\text{i.i.d.}}{\sim} N$, where $L$ and $N$ are positive discrete random variables, with $P(L > 1) > 0$ and $P(N > 1) > 0$.

The observation times and measurements are assumed to be independent of the number of measurements, that is, for any subsets $\mathcal{L}_i \subseteq \{1, \ldots, L_i\}$ and $\mathcal{J}_i \subseteq \{1, \ldots, N_i\}$, and for all $i = 1, \ldots, n$,

(B1.2) $(\{S_{il} : l \in \mathcal{L}_i\}, \{U_{il} : l \in \mathcal{L}_i\})$ is independent of $L_i$, and $(\{T_{ij} : j \in \mathcal{J}_i\}, \{Y_{ij} : j \in \mathcal{J}_i\})$ is independent of $N_i$.

Let $K_1(\cdot)$ and $K_2(\cdot, \cdot)$ be nonnegative univariate and bivariate kernel functions that are used in the smoothing steps for the mean functions $\mu_X$, $\mu_Y$, covariance surfaces $G_X$, $G_Y$, and cross-covariance structure $C$. Assume that $K_1$ and $K_2$ are compactly supported densities with zero means and finite variances. Let $b_X = b_X(n)$, $b_Y = b_Y(n)$, $h_X = h_X(n)$, $h_Y = h_Y(n)$ be the bandwidths for estimating $\hat{\mu}_X$ and $\hat{\mu}_Y$ (30), $\hat{G}_X$ and $\hat{G}_Y$ (31), and $h_1 = h_1(n)$, $h_2 = h_2(n)$ be the bandwidths for obtaining $\hat{C}$ (36). We develop asymptotics as the number of subjects $n \to \infty$, and require the following:

(B2.1) $b_X \to 0$, $b_Y \to 0$, $nb_X^4 \to \infty$, $nb_Y^4 \to \infty$, $nb_X^6 < \infty$ and $nb_Y^6 < \infty$.

(B2.2) $h_X \to 0$, $h_Y \to 0$, $nh_X^6 \to \infty$, $nh_Y^6 \to \infty$, $nh_X^8 < \infty$ and $nh_Y^8 < \infty$.

(B2.3) Without loss of generality, $h_1/h_2 \to 1$, $nh_1^6 \to \infty$ and $nh_1^8 < \infty$.

Define the Fourier transformations of $K_1(u)$, $K_2(u, v)$ by $\kappa_1(t) = \int e^{-iut} K_1(u) \, du$ and $\kappa_2(t, s) = \int e^{-(iut+ivs)} K_2(u, v) \, du \, dv$. They satisfy the following:

(B3.1) $\kappa_1(t)$ is absolutely integrable, that is, $\int |\kappa_1(t)| \, dt < \infty$.

(B3.2) $\kappa_2(t, s)$ is absolutely integrable, that is, $\int \int |\kappa_2(t, s)| \, dt \, ds < \infty$.

Assume that the fourth moments of $Y$ and $U$, centered at $\mu_Y(T)$ and $\mu_X(S)$, are finite, that is,

(B4) $E[(Y - \mu_Y(T))^4] < \infty$, $E[(U - \mu_X(S))^4] < \infty$.

Let $S_1$ and $S_2$ be i.i.d. as $S$, and $U_1$ and $U_2$ be the repeated measurements of $X$ made on the same subject, taken at $S_1$ and $S_2$ separately. Assume $(S_{il_1}, S_{il_2}, U_{il_1}, U_{il_2})$, $1 \leq l_1 \neq l_2 \leq L_i$, is identically distributed as $(S_1, S_2, U_1, U_2)$ with joint density function $g_X(s_1, s_2, u_1, u_2)$, and analogously for $(T_{ij_1}, T_{ij_2}, V_{ij_1}, V_{ij_2})$ with identical joint density function $g_Y(t_1, t_2, v_1, v_2)$. Appropriate regularity assumptions are imposed for the marginal and joint densities, $f_S(s)$, $f_T(t)$, $g_1(s, u)$, $g_2(t, v)$, $g_X(s_1, s_2, u_1, u_2)$ and $g_Y(t_1, t_2, v_1, v_2)$.

Define the rank one operator $f \otimes g = \langle f, h \rangle y$, for $f, h \in H$, and denote the separable Hilbert space of Hilbert–Schmidt operators on $H$ by $F \equiv \sigma_2(H)$, endowed by $\langle T_1, T_2 \rangle_F = \operatorname{tr}(T_1 T_2^*) = \sum_j \langle T_1 u_j, T_2 u_j \rangle_H$ and $\|T\|_F^2 = \langle T, T \rangle_F$, where $T_1, T_2, T \in F$, and $\{u_j : j \geq 1\}$ is any complete orthonormal system in $H$. The covariance operator $\mathbf{G}_X$ (resp. $\widehat{\mathbf{G}}_X$) is generated by the kernel $G_X$ (resp. $\widehat{G}_X$), that is, $\mathbf{G}_X(f) = \int_{\mathcal{S}} G_X(s, t) f(s) \, ds$, $\widehat{\mathbf{G}}(f) = \int_{\mathcal{S}} \widehat{G}_X(s, t) f(s) \, ds$. Then $\mathbf{G}_X$ and $\widehat{\mathbf{G}}_X$ are Hilbert–Schmidt operators, and Theorem 1 in [31] implies that $\|\widehat{\mathbf{G}}_X - \mathbf{G}_X\|_F = O_p(1/(\sqrt{n} h_X^2))$.

Let $\mathcal{I}_i = \{j : \rho_j = \rho_i\}$, $\mathcal{I}' = \{i : |\mathcal{I}_i| = 1\}$, where $|\mathcal{I}_i|$ denotes the number of elements in $\mathcal{I}_i$. Let $\mathbf{P}_j^X = \sum_{k \in \mathcal{I}_j} \psi_k \otimes \psi_k$, and $\widehat{\mathbf{P}}_j^X = \sum_{m \in \mathcal{I}_j} \hat{\psi}_m \otimes \hat{\psi}_m$ denote the true and estimated orthogonal projection operators from $H$ to the subspace spanned by $\{\psi_m : m \in \mathcal{I}_j\}$. For fixed $j$, let

$$\delta_j^X = \tfrac{1}{2} \min\{|\rho_l - \rho_j| : l \notin \mathcal{I}_j\},$$

and let $\mathbf{\Lambda}_{\delta_j^X} = \{z \in \mathcal{C} : |z - \rho_j| = \delta_j^X\}$, where $\mathcal{C}$ stands for the set of complex numbers. The resolvent of $\mathbf{G}_X$ (resp. $\widehat{\mathbf{G}}_X$) is denoted by $\mathbf{R}_X$ (resp. $\widehat{\mathbf{R}}_X$), that is, $\mathbf{R}_X(z) = (\mathbf{G}_X - zI)^{-1}$ [resp. $\widehat{\mathbf{R}}_X(z) = (\widehat{\mathbf{G}}_X - zI)^{-1}$]. Let

$$A_{\delta_j^X} = \sup\{\|\mathbf{R}_X(z)\|_F : z \in \mathbf{\Lambda}_{\delta_j^X}\},$$

and analogously define sequences $\{\delta_j^Y\}$ and $\{A_{\delta_j^Y}\}$ for the response process $Y$. We assume that the numbers $M = M(n)$ and $K = M(n)$ of included eigenfunctions depend on the sample size $n$, such that as $n \to \infty$, if $M(n) \to \infty$ and $K(n) \to \infty$,

$$
\begin{align}
\tau_n &= \sum_{m=1}^{M(n)} \frac{\delta_m^X A_{\delta_m^X}}{\sqrt{n} h_X^2 - A_{\delta_m^X}} \to 0, \\
\upsilon_n &= \sum_{k=1}^{K(n)} \frac{\delta_k^Y A_{\delta_k^Y}}{\sqrt{n} h_Y^2 - A_{\delta_k^Y}} \to 0, \\
\varsigma_n &= \frac{KM}{\sqrt{n} h_1 h_2} \to 0.
\end{align}
$$

(B5)

The main effect of this assumption is to impose certain constraints on the rate of $K$ and $M$ in relation to $n$ and the bandwidths.

We denote the remainder in (A7) by

$$(29) \qquad \pi_n = \sup_{t \in \mathcal{T}} \left| \sum_{m=M(n)}^{\infty} \sum_{k,\ell=K(n)}^{\infty} \frac{\sigma_{km} \sigma_{\ell m}}{\rho_m} \phi_k(t) \phi_\ell(t) \right|.$$

**A.2. Estimation procedures.** We begin by summarizing the estimation procedures for the components of models (1) and (2) in the following; see Yao, Müller and Wang [31] for further details. Define the local linear scatterplot smoothers for $\mu_X(s)$ through minimizing

$$(30) \qquad \sum_{i=1}^{n} \sum_{l=1}^{L_i} K_1 \left( \frac{S_{il} - s}{b_X} \right) \{ U_{ij} - \beta_0^X - \beta_1^X (s - S_{il}) \}^2,$$

with respect to $\beta_0^X$, $\beta_1^X$, leading to $\hat{\mu}_X(s) = \hat{\beta}_0^X(s)$.

Let $G_i^X(S_{il_1}, S_{il_2}) = (U_{il_1} - \hat{\mu}_X(S_{il_1}))(U_{il_2} - \hat{\mu}_X(S_{il_2}))$, and define the local linear surface smoother for $G_X(s_1, s_2)$ through minimizing

$$(31) \qquad \sum_{i=1}^{n} \sum_{1 \le l_1 \ne l_2 \le L_i} K_2 \left( \frac{S_{il_1} - s_1}{h_X}, \frac{S_{il_2} - s_2}{h_X} \right)$$
$$\times \{ G_i^X(S_{il_1}, S_{il_2}) - f(\beta^X, (s_1, s_2), (S_{il_1}, S_{il_2})) \}^2,$$

where $f(\beta^X, (s_1, s_2), (S_{il_1}, S_{il_2})) = \beta_0^X + \beta_{11}^X(s_1 - S_{il_1}) + \beta_{12}^X(s_2 - S_{il_2})$, with respect to $\beta^X = (\beta_0^X, \beta_{11}^X, \beta_{12}^X)$, yielding $\hat{G}_X(s_1, s_2) = \hat{\beta}_0^X(s_1, s_2)$.

For the estimation of $\sigma_X^2$, we fit a local quadratic component orthogonal to the diagonal of $G_X$, and a local linear component in the direction of the diagonal. Denote the diagonal of the resulting surface estimate by $\tilde{G}_X(s)$, and a local linear smoother focusing on diagonal values $\{ G_X(s,s) + \sigma_X^2 \}$ by $\hat{V}_X(s)$. Let $a_X = \inf \{ s : s \in \mathcal{S} \}$, $b_X = \sup \{ s : s \in \mathcal{S} \}$, $|\mathcal{S}| = b_X - a_X$, $\mathcal{S}_1 = [a_X + |\mathcal{S}|/4, b_X - |\mathcal{S}|/4]$. The estimate of $\sigma_X^2$ is

$$(32) \qquad \hat{\sigma}_X^2 = 2 \int_{\mathcal{S}_1} \{ \hat{V}_X(s) - \tilde{G}_X(s) \} \, ds / |\mathcal{S}|,$$

if $\hat{\sigma}_X^2 > 0$ and $\hat{\sigma}_X^2 = 0$ otherwise.

The estimates of $\{ \rho_m, \psi_m \}_{m \ge 1}$ correspond to the solutions $\{ \hat{\rho}_m, \hat{\psi}_m \}_{m \ge 1}$ of the eigenequations

$$(33) \qquad \int_{\mathcal{S}} \hat{G}_X(s_1, s_2) \hat{\psi}_m(s_1) \, ds_1 = \hat{\rho}_m \hat{\psi}_m(s_2),$$

with orthonormal constraints on $\{ \hat{\psi}_m \}_{m \ge 1}$.

Let $\hat{\mu}_X^{(-i)}$ and $\hat{\psi}_m^{(-i)}$ be the estimated mean and eigenfunctions after removing the data for $X_i$. One-curve-leave-out cross-validation aims to minimize

$$(34) \qquad \mathrm{CV}_X(M) = \sum_{i=1}^{n} \sum_{l=1}^{L_i} \{ U_{il} - \hat{X}_i^{(-i)}(S_{il}) \}^2$$

with respect to $M$, where $\widehat{X}_i^{(-i)}(s) = \hat{\mu}_X^{(-i)}(s) + \sum_{m=1}^{M} \hat{\zeta}_{im}^{(-i)} \hat{\psi}_m^{(-i)}(s)$, and $\hat{\zeta}_{im}^{(-i)}$ is obtained by (11). The AIC criterion as a function of $M$ is given by

$$
\begin{aligned}
\text{AIC}(M) = \sum_{i=1}^{n} \Bigg\{ \frac{1}{2\hat{\sigma}_X^2} \Bigg( \widetilde{U}_i - \hat{\mu}_{X_i} - \sum_{m=1}^{M} \hat{\zeta}_{im} \hat{\psi}_{im} \Bigg)^T \\
\times \Bigg( \widetilde{U}_i - \hat{\mu}_{X_i} - \sum_{m=1}^{M} \hat{\zeta}_{im} \hat{\psi}_{im} \Bigg) \\
+ \frac{L_i}{2} \log(2\pi) + \frac{L_i}{2} \log \hat{\sigma}_X^2 \Bigg\} + M,
\end{aligned}
\tag{35}
$$

where $\widetilde{U}_i = (U_{i1}, \ldots, U_{iL_i})^T$, $\hat{\mu}_{X_i} = (\hat{\mu}_X(S_{i1}), \ldots, \hat{\mu}_X(S_{iL_i}))^T$, $\hat{\psi}_{im} = (\hat{\psi}_m(S_{i1}), \ldots, \hat{\psi}_m(S_{iL_i}))^T$, and $\hat{\zeta}_{im}$ is obtained by (11). We proceed analogously for the corresponding estimates for the components of model (2) regarding the response process $Y$.

The local linear surface smoother for the cross-covariance surface $C(s,t)$ is obtained through minimizing

$$
\sum_{i=1}^{n} \sum_{j=1}^{N_i} \sum_{l=1}^{L_i} K_2 \bigg( \frac{S_{il} - s}{h_1}, \frac{T_{ij} - t}{h_2} \bigg) \{ C_i(S_{il}, T_{ij}) - f(\beta, (s,t), (S_{il}, T_{ij})) \}^2
\tag{36}
$$

with regard to $\beta = (\beta_0, \beta_{11}, \beta_{12})$, leading to $\widehat{C}(s,t) = \hat{\beta}_0(s,t)$.

**A.3. Preliminary consistency results.** Applying Theorems 1 and 2 of [31], we obtain uniform consistency of the estimate of the mean functions, covariance functions, eigenvalues and eigenfunctions of the predictor and response processes. Under assumption (A2.3), this extends to the cross-covariance function.

LEMMA A.1. *Under* (B1.1)–(B5), *and appropriate regularity assumptions for* $f_S(s)$, $f_T(t)$, $g_1(s,u)$, $g_2(t,v)$, $g_X(s_1, s_2, u_1, u_2)$ *and* $g_Y(t_1, t_2, v_1, v_2)$,

$$
|\hat{\rho}_m - \rho_m| = O_p \bigg( \frac{\delta_m^X A_{\delta_m^X}}{\sqrt{n} h_X^2 - A_{\delta_m^X}} \bigg), \qquad |\hat{\lambda}_k - \lambda_k| = O_p \bigg( \frac{\delta_k^Y A_{\delta_k^Y}}{\sqrt{n} h_Y^2 - A_{\delta_k^Y}} \bigg).
\tag{37}
$$

*Considering eigenvalues* $\rho_m$ *and* $\lambda_k$ *of multiplicity one respectively,* $\hat{\psi}_m$ *and* $\hat{\phi}_k$ *can be chosen such that*

$$
\begin{aligned}
\sup_{s \in \mathcal{S}} |\hat{\psi}_m(s) - \psi_m(s)| = O_p \bigg( \frac{\delta_m^X A_{\delta_m^X}}{\sqrt{n} h_X^2 - A_{\delta_m^X}} \bigg), \\
\sup_{t \in \mathcal{T}} |\hat{\phi}_k(t) - \phi_k(t)| = O_p \bigg( \frac{\delta_k^Y A_{\delta_k^Y}}{\sqrt{n} h_Y^2 - A_{\delta_k^Y}} \bigg),
\end{aligned}
\tag{38}
$$

*and, furthermore,*

$$\text{(39)} \qquad \sup_{(s,t)\in\mathcal{S}\times\mathcal{T}} |\widehat{C}(s,t) - C(s,t)| = O_p\bigg(\frac{1}{\sqrt{n}h_1 h_2}\bigg).$$

*As a consequence of* (38) *and* (39),

$$\text{(40)} \quad |\hat{\sigma}_{km} - \sigma_{km}| = O_p\bigg(\max\bigg\{\frac{\delta_m^X A_{\delta_m^X}}{\sqrt{n}h_X^2 - A_{\delta_m^X}}, \frac{\delta_k^Y A_{\delta_k^Y}}{\sqrt{n}h_Y^2 - A_{\delta_k^Y}}, \frac{1}{\sqrt{n}h_1 h_2}\bigg\}\bigg),$$

*where the* $O_p(\cdot)$ *terms in* (37), (38) *and* (40) *hold uniformly over all* $k$ *and* $m$.

PROOF. Part of the proof follows that of Theorem 2 in [31]. Additional arguments are needed to obtain the convergence rates in (38). Since $\widehat{\mathbf{R}}_X(z) = \mathbf{R}_X(z)[I + (\widehat{\mathbf{G}}_X - \mathbf{G}_X)\mathbf{R}_X(z)]^{-1} = \mathbf{R}_X(z)\sum_{l=0}^{\infty}[(\widehat{\mathbf{G}}_X - \mathbf{G}_X)\mathbf{R}_X(z)]^l$, $\|\widehat{\mathbf{R}}_X(z) - \mathbf{R}_X(z)\|_F \leq (\|\widehat{\mathbf{G}}_X - \mathbf{G}_X\|_F\|\mathbf{R}_X(z)\|_F)/(1 - \|\widehat{\mathbf{G}}_X - \mathbf{G}_X\|_F\|\mathbf{R}_X(z)\|_F)$. Note that $\mathbf{P}_j^X = (2\pi i)^{-1}\int_{\mathbf{\Lambda}_{\delta_j^X}}\mathbf{R}_X(z)\,dz$, $\widehat{\mathbf{P}}_j^X = (2\pi i)^{-1}\int_{\mathbf{\Lambda}_{\delta_j^X}}\widehat{\mathbf{R}}_X(z)\,dz$. Therefore,

$$\|\widehat{\mathbf{P}}_j^X - \mathbf{P}_j^X\|_F \leq \int_{\mathbf{\Lambda}_{\delta_j^X}} \|\widehat{\mathbf{R}}_X(z) - \mathbf{R}_X(z)\|_F\,dz/(2\pi)$$

$$\leq \delta_j^X \frac{\|\widehat{\mathbf{G}}_X - \mathbf{G}_X\|_F A_{\delta_j^X}}{1 - \|\widehat{\mathbf{G}}_X - \mathbf{G}_X\|_F A_{\delta_j^X}}$$

$$\leq \frac{\delta_j^X A_{\delta_j^X}}{\sqrt{n}h_X^2 - A_{\delta_j^X}}.$$

Applying the arguments used in the proof of Theorem 2 in [31] leads to (37) and (38). Under (A2.3), the uniform convergence of $\widehat{C}$ is an extension of Theorem 1 in [31]. Then (40) follows by applying (38) and (39). □

LEMMA A.2. *Under* (A1) *the right-hand side of* (5) *converges in the* $L^2$ *sense. Furthermore, under* (A2) *the right-hand side of* (5) *converges uniformly on* $\mathcal{S}\times\mathcal{T}$.

PROOF. Let $\beta_{KM}(s,t) = \sum_{k=1}^{K}\sum_{m=1}^{M}\sigma_{km}\psi_m(s)\phi_k(t)/\rho_m$. Observing the orthonormality of $\{\psi_m\}_{m\geq 1}$ and $\{\phi_k\}_{k\geq 1}$, $\int_{\mathcal{T}}\int_{\mathcal{S}}\beta_{KM}^2(s,t)\,dt = \sum_{k=1}^{K}\sum_{m=1}^{M}\sigma_{km}^2/\rho_m^2$, and it is obvious that $\beta_{KM}$ converges in the $L^2$ sense under (A1), that is, $\int_{\mathcal{T}}\int_{\mathcal{S}}[\beta_{KM}(s,t) - \beta(s,t)]^2\,ds\,dt \to 0$ as $K, M \to \infty$. Let $\beta = \lim_{K,M\to\infty}\beta_{KM}$. For all $(s,t)\in\mathcal{S}\times\mathcal{T}$, let $\gamma_{KM}(s,t) = \sum_{k=1}^{K}\sum_{m=1}^{M}|\sigma_{km}\psi_m(s)\phi_k(t)|/\rho_m$. By (A2), the monotonically increasing net of continuous real-valued functions $\{\gamma_{KM}(s,t)\}$ converges pointwise to a continuous function $\gamma(s,t)$. By applying Dini's theorem, the net $\{\gamma_{KM}(s,t)\}$ converges to $\gamma(s,t)$ uniformly on the compact set $\mathcal{S} \times \mathcal{T}$, which implies the uniform convergence of the right-hand side of (5).

We denote the remainders in (A1) and (A2) as $M(n), K(n) \to \infty$ as

$$\theta_n = \sum_{k=K(n)}^{\infty} \sum_{m=M(n)}^{\infty} \frac{\sigma_{km}^2}{\rho_m^2},$$

(41)

$$\vartheta_n = \sup_{s,t \in \mathcal{T} \times \mathcal{S}} \left| \sum_{k=K(n)}^{\infty} \sum_{m=M(n)}^{\infty} \sigma_{km} \psi_m(s) \phi_k(t)/\rho_m \right|. \qquad \square$$

LEMMA A.3. *If* (A4) *holds, as* $K, M \to \infty$,

$$\sup_{t \in \mathcal{T}} E[Y_{KM}^*(t) - E[Y^*(t)|X^*]]^2 \longrightarrow 0, \tag{42}$$

*and for given* $L^*$ *and* $S^*$,

$$\sup_{t \in \mathcal{T}} E[\widetilde{Y}_{KM}^*(t) - \widetilde{Y}^*(t)]^2 \longrightarrow 0. \tag{43}$$

PROOF. To prove (42), note that

$$\sup_{t \in \mathcal{T}} E[Y_{KM}^*(t) - E[Y^*(t)|X^*]]^2 = \sum_{m=M+1}^{\infty} \sum_{k=K+1}^{\infty} \frac{\sigma_{km}^2}{\rho_m \lambda_k} \sup_{t \in \mathcal{T}} \{\lambda_k \phi_k^2(t)\}.$$

From the Karhunen–Loève theorem, we know that $\sum_{k=1}^{\infty} \lambda_k \phi(s)\phi(t)$ converges uniformly in $s, t \in \mathcal{T}$, which implies that $\sup_{t \in \mathcal{T}} \lambda_k \phi_k(t)^2$ converges to zero as $k \to \infty$. If (A4) holds, then (42) follows. Given $L^*$ and $S^*$, since $\widetilde{Y}_{KM}^*(t) - \widetilde{Y}^*(t) = E[\sum_{m=M+1}^{\infty} \sum_{k=K+1}^{\infty} \sigma_{km} \phi_k(t) \zeta_m^*/\rho_m | \widetilde{U}^*]$,

$$\sup_{t \in \mathcal{T}} \left\{ E\left[ E\left[ \sum_{M+1}^{\infty} \sum_{K+1}^{\infty} \frac{\sigma_{km} \phi_k(t) \zeta_m^*}{\rho_m} \bigg| \widetilde{U}^* \right]^2 \right] + E\left[ \operatorname{var}\left( \sum_{M+1}^{\infty} \sum_{K+1}^{\infty} \frac{\sigma_{km} \phi_k(t) \zeta_m^*}{\rho_m} \bigg| \widetilde{U}^* \right) \right] \right\}$$

$$= \sum_{M+1}^{\infty} \sum_{K+1}^{\infty} \frac{\sigma_{km}^2}{\rho_m \lambda_k} \sup_{t \in \mathcal{T}} \{\lambda_k \phi_k^2(t)\} \longrightarrow 0,$$

and

$$E\left[ \operatorname{var}\left( \sum_{M+1}^{\infty} \sum_{K+1}^{\infty} \sigma_{km} \phi_k(t) \zeta_m^*/\rho_m | \widetilde{U}^* \right) \right] \geq 0 \qquad \text{for all } t \in \mathcal{T}.$$

The result (43) follows.  $\square$

### A.4. Proofs of the main results.

PROOF OF THEOREM 1. To prove (22), observing the orthonormality of the eigenfunction basis,

$$\int_{\mathcal{T}}\int_{\mathcal{S}}[\hat{\beta}(s,t)-\beta(s,t)]^2\,ds\,dt$$

$$=\int_{\mathcal{T}}\int_{\mathcal{S}}\left\{\sum_{k=1}^{K-1}\sum_{m=1}^{M-1}\left[\frac{\hat{\sigma}_{km}}{\hat{\rho}_m}\hat{\psi}_m(s)\hat{\phi}_k(t)-\frac{\sigma_{km}}{\rho_m}\psi_m(s)\phi_k(t)\right]\right\}^2\,ds\,dt$$

$$+\sum_{k=K}^{\infty}\sum_{m=M}^{\infty}\frac{\sigma_{km}^2}{\rho_m^2}$$

$$+\int_{\mathcal{T}}\int_{\mathcal{S}}\left[\sum_{k=K}^{\infty}\sum_{m=M}^{\infty}\frac{\sigma_{km}}{\rho_m}\psi_m(s)\phi_k(t)\right]$$

$$\times\left\{\sum_{k=1}^{K-1}\sum_{m=1}^{M-1}\left[\frac{\hat{\sigma}_{km}}{\hat{\rho}_m}\hat{\psi}_m(s)\hat{\phi}_k(t)-\frac{\sigma_{km}}{\rho_m}\psi_m(s)\phi_k(t)\right]\right\}\,ds\,dt$$

$$\equiv Q_1(n)+Q_2(n)+Q_3(n).$$

Then (A1) and (B5) imply that $Q_2(n)\to 0$ as $M(n),K(n)\to\infty$, that is, $n\to\infty$. By (37), (38), (40) and (B5), Slutsky's theorem implies

$$Q_1(n)=O_p\left(\sum_{m=1}^{M}\frac{\delta_m^X A_{\delta_m^X}}{\sqrt{n}h_X^2-A_{\delta_m^X}}+\sum_{k=1}^{K}\frac{\delta_k^Y A_{\delta_k^Y}}{\sqrt{n}h_Y^2-A_{\delta_k^Y}}+\frac{KM}{\sqrt{n}h_1h_2}\right)\xrightarrow{p}0,$$

(44)

as $n\to\infty$. For $Q_3(n)$, using the Cauchy–Schwarz inequality, $Q_3(n)^2\le Q_1(n)\times Q_2(n)$, thus $Q_3(n)\xrightarrow{p}0$ as $n\to\infty$. Then (22) follows. To prove (23), note that

$$\sup_{s,t}|\hat{\beta}(s,t)-\beta(s,t)|\le\sup_{s,t}\left|\sum_{k=1}^{K-1}\sum_{m=1}^{M-1}\left[\frac{\hat{\sigma}_{km}}{\hat{\rho}_m}\hat{\psi}_m(s)\hat{\phi}_k(t)-\frac{\sigma_{km}}{\rho_m}\psi_m(s)\phi_k(t)\right]\right|$$

$$+\sup_{s,t}\left|\sum_{k=K}^{\infty}\sum_{m=M}^{\infty}\frac{\sigma_{km}}{\rho_m}\psi_m(s)\phi_k(t)\right|$$

$$\equiv Q_4(n)+Q_5(n).$$

One has $Q_5(n)\to 0$ as $M(n),K(n)\to\infty$ if (A2) holds. Observing (37), (38), (40) and (B5), $Q_4(n)\xrightarrow{p}0$ as $n\to\infty$, leading to (23). □

PROOF OF THEOREM 2. For given $L^*$ and $S^*$, define $\widetilde{Y}_{KM}^*(t)=\mu_Y(t)+\sum_{k=1}^{K}\sum_{m=1}^{M}\sigma_{km}\tilde{\zeta}_m^*\phi_k(t)/\rho_m$, where $\tilde{\zeta}_m^*$ is defined in (10). Recall $\widetilde{Y}^*(t) = \mu_Y(t) + \sum_{k=1}^{\infty}\sum_{m=1}^{\infty}\sigma_{km}\tilde{\zeta}^*_m\phi_k(t)/\rho_m$, $\widehat{Y}^*_{KM}(t) = \hat{\mu}_Y(t) + \sum_{k=1}^{K}\sum_{m=1}^{M}\hat{\sigma}_{km}\hat{\zeta}^*_m\hat{\phi}_k(t)/\hat{\rho}_m$, with $\hat{\zeta}^*_m$ defined in (11). Note that

$$|\widehat{Y}^*_{KM}(t) - \widetilde{Y}^*(t)| \leq |\widehat{Y}^*_{KM}(t) - \widetilde{Y}^*_{KM}(t)| + |\widetilde{Y}^*_{KM}(t) - \widetilde{Y}^*(t)|.$$

Lemma A.3 implies $\widetilde{Y}^*_{KM}(t) \overset{p}{\to} \widetilde{Y}^*(t)$ as $K, M \to \infty$ and $n \to \infty$. Applying Theorem 1 in [31], Lemma A.1 and (B5), one has $\sup_{t\in\mathcal{T}}|\hat{\mu}_Y(t) - \mu_Y(t)| = O_p(1/(\sqrt{n}b_Y))$ and $|\hat{\zeta}^*_m - \tilde{\zeta}^*_m| = O_p(\delta^X_m A_{\delta^X_m}/(\sqrt{n}h^2_X - A_{\delta^X_m}))$ as $n \to \infty$. Then $\sup_{t\in\mathcal{T}}|\widehat{Y}^*_{KM}(t) - \widetilde{Y}^*_{KM}(t)| \overset{p}{\to} 0$ as $n \to \infty$ by Slutsky's theorem, and (24) follows.  □

PROOF OF THEOREM 3.  Analogous to the proof of Theorem 4 in [31] with slight modifications similar to the arguments used in the proof of Theorem 1.  □

PROOF OF THEOREM 4.  To prove (26), note

$$\left|\frac{\sum_{m=1}^{M}\sum_{k=1}^{K}\hat{\sigma}^2_{km}/\hat{\rho}}{\sum_{k=1}^{K}\hat{\lambda}_k} - \frac{\sum_{m,k=1}^{\infty}\sigma^2_{km}/\rho_m}{\sum_{k=1}^{\infty}\lambda_k}\right|$$

$$\leq \left|\frac{\sum_{k=1}^{K}\sum_{m=1}^{M}\hat{\sigma}^2_{km}/\hat{\rho}_m}{\sum_{k=1}^{K}\hat{\lambda}_k} - \frac{\sum_{k=1}^{K}\sum_{m=1}^{M}\sigma^2_{km}/\rho_m}{\sum_{k=1}^{K}\lambda_k}\right|$$

$$+ \left|\frac{\sum_{k=1}^{K}\sum_{m=1}^{M}\sigma^2_{km}/\rho_m}{\sum_{k=1}^{K}\lambda_k} - \frac{\sum_{m,k=1}^{\infty}\sigma^2_{km}/\rho_m}{\sum_{k=1}^{\infty}\lambda_k}\right|$$

$$\equiv Q_1(n) + Q_2(n).$$

Since the nonnegative series $\sum_{m,k=1}^{\infty}\sigma^2_{km}/\rho_m \leq \sum_{k=1}^{\infty}\lambda_k < \infty$, we have $Q_2(n) \to 0$ as $M, K \to \infty$. Observing (37), (38), (40) and (B5), one finds that $Q_1(n) \overset{p}{\to} 0$ as $n \to \infty$, by applying similar arguments to (44), leading to (26). It is obvious that (27) implies (28). To show (27), let $\nu_{KM}(t) = \sum_{m=1}^{M}\sum_{k,\ell=1}^{K}|\sigma_{km}\sigma_{\ell m}\phi_k(t)\phi_\ell(t)|/\rho_m$ for all $t \in \mathcal{T}$. By (A7) and Dini's theorem, the net $\{\nu_{KM}(t)\}$ converges to $\nu(t)$ uniformly on the compact set $\mathcal{T}$, which implies the uniform convergence of $R^2_{KM}(t)$ as $K, M \to \infty$. Applying arguments similar to those used to prove (26), one obtains (27).  □

**Acknowledgments.**  We wish to thank two referees and an Associate Editor for their most helpful and insightful remarks, and are indebted to Dr. C. H. Morrell, Loyola College, for his support and help to access an anonymous version of a subset of the BLSA data.

F. Yao
Department of Statistics
Colorado State University
Fort Collins, Colorado 80523
USA

H.-G. Müller
J.-L. Wang
Department of Statistics
University of California, Davis
Davis, California 95616
USA
e-mail: mueller@wald.ucdavis.edu